\documentclass[letterpaper, paper, 10pt]{AAS}	

\usepackage{bm}
\usepackage{amsmath}
\usepackage{siunitx} 
\usepackage[colorlinks=true, pdfstartview=FitV, linkcolor=black, citecolor= black, urlcolor= black]{hyperref}
\usepackage{overcite}
\usepackage{footnpag}			      	

\usepackage{gensymb}
\DeclareMathAlphabet{\mathpzc}{OT1}{pzc}{m}{it}
\usepackage{mathtools}
\usepackage{adjustbox}
\usepackage{array}
\usepackage{graphicx}
\usepackage{caption}
\usepackage{hyperref}
\usepackage{amssymb}
\usepackage[super]{nth}
\usepackage{color}
\usepackage{soul}
\usepackage[bottom]{footmisc}
\usepackage{tabularx}
\usepackage[table,xcdraw]{xcolor}
\usepackage{textcomp}
\usepackage{dcolumn}
\usepackage{fancyhdr}
\usepackage{threeparttable}
\usepackage{xspace}
\usepackage{subfigure}
\usepackage{pictex}
\usepackage{color}
\usepackage{soul}
\usepackage{siunitx}
\usepackage{gensymb}
\usepackage{multirow}

\usepackage{caption}
\usepackage{setspace}
\usepackage{algorithm}
\usepackage{algpseudocode}

\def\bfr{\bm{r}}
\def\bfv{\bm{v}}
\def\bfp{\bm{p}}
\def\bft{\bm{t}}
\def\bftau{\bm{\tau}}
\def\bfPi{\bm{\Pi}}

\newcommand{\norm}[1]{\left\|#1\right\|}

\DeclareMathOperator{\acos}{acos}

\PaperNumber{19-701}

\title{\textbf{
A Time-Dependent TSP Formulation for the Design of an Active Debris Removal Mission using Simulated Annealing 
}}
\author{Lorenzo Federici\footnote{PhD student, Department of Mechanical and Aerospace Engineering, Sapienza University of Rome, Via Eudossiana 18 - 00184, Rome, Italy. lorenzo.federici@uniroma1.it}, Alessandro Zavoli\footnote{Research Assistant, Department of Mechanical and Aerospace Engineering, Sapienza University of Rome, Via Eudossiana 18 - 00184, Rome, Italy. alessandro.zavoli@uniroma1.it}, Guido Colasurdo\footnote{Full Professor, Department of Mechanical and Aerospace Engineering, Sapienza University of Rome, Via Eudossiana 18 - 00184, Rome, Italy. guido.colasurdo@uniroma1.it}}

\begin{document}

\maketitle

\begin{abstract}
This paper proposes a formulation of the Active Debris Removal (ADR) Mission Design problem as a modified Time-Dependent Traveling Salesman Problem (TDTSP). The TDTSP is a well-known combinatorial optimization problem, whose solution is the cheapest mono-cyclic tour connecting a number of non-stationary cities in a map. The problem is tackled with an optimization procedure based on Simulated Annealing, that efficiently exploits a natural encoding and a careful choice of mutation operators. The developed algorithm is used to simultaneously optimize the targets sequence and the rendezvous epochs of an impulsive ADR mission. Numerical results are presented for sets comprising up to 20 targets.
\end{abstract}

\section{Introduction}


The Traveling Salesman Problem (TSP) is a well-known combinatorial optimization problem, whose solution is the cheapest tour which allows a salesman to visit, only once, a number of cities in a map; the cost of each city-to-city transfer is, typically, the traveled distance or the fuel consumption. So, given a set of $n$ cities, with cost $c_{ij}$ associated to the leg from city $i$ to $j$, the optimal tour is represented by the permutation $\pi$ of the numbers $\{1,\dots,n\}$ which minimizes $\sum_{i = 1}^{n - 1} c_{\pi(i) \pi(i+1)} $. 
In the time-dependent version of the problem (TDTSP), the cost associated with traveling  between any two cities changes with time. 
Among those present in literature, two TDTSP variants, with increasing complexity, are here considered: in approach (A), the cost $c_{ij}$ to go from city $i$ to city $j$ is dependent only on the position of the leg $i$-$j$ in the whole tour; in approach (B), the cost $c_{ij}$ depends both on the starting time and on the duration of the transfer $i$-$j$.

Active Debris Removal (ADR) missions can be seen as peculiar instances of the TDTSP, 
where  an active (chaser) spacecraft is asked to visit, that is, to perform a rendezvous, with a certain number of targets (space debris), making the best use of the on-board propellant.
Such kind of missions are increasing in popularity among space agencies all over the world, as the sustainability of the extra-atmospheric environment is becoming compromised by the huge amount of ``space garbage'' now orbiting Earth.\cite{liou2008instability}
In addition, some
companies and organizations plan to launch mega-constellations of telecommunications satellites
over the next decade\cite{foreman2017large}, so the need for space debris mitigation efforts is becoming increasingly urgent. 
A cost-competitive space program would involve the removal of several dozens of small debris with each single mission; such a complex scenario could became feasible only with the best possible use of the propellant on-board of the chaser spacecraft. 
As a consequence, a well-designed ADR mission would require the optimization of a multi-target
rendezvous trajectory.

This multi-target rendezvous problem can be seen as an extension of the
single-target time-fixed rendezvous, a well-known problem in spaceflight mechanics,
that have been extensively dealt in literature assuming either a finite \cite{Lu2012} or an impulsive thrust \cite{Prussing1986} model. In the latter case, four burns permit the achievement of the optimal solution, if close coplanar orbits are considered \cite{prussing1967optimal}.
However, the combinatorial ``flavor'' introduced by the possibility to reorder the target sequence makes the multi-target rendezvous problem much closer to a TDTSP instance than to a traditional space trajectory optimization problem, as the choice of the sequence of targets may affect the overall propellant consumption more than the use of a sub-optimal flight strategy.

A number of authors dealt with long term or time-free ADR missions  aimed at removing a small number of debris from Sun synchronous orbits (at a rate of three to ten per year). 
These missions heavily rely on  $J_2$ orbital perturbation for the alignment of the orbital planes of consecutive targets before starting  the rendezvous maneuver, in order to reduce the mission cost.\cite{ shen2017optimization} 
However, such an operational scenario may become impractical in presence of strict time-constraints or long debris sequences.
Some researchers focused on 
selecting the trajectory that maximizes the number of removed debris within a very large catalogue, while complying with some propellant constraints.\cite{Barea2009}
However, debris features such as size, orbit, or proximity to relevant active spacecrafts may considerably increase the importance of removing specific debris, rather than the maximum number of debris in a huge set.
In this respect,
this paper focuses on the design of an ADR mission involving a prescribed set of debris, that move on the same orbital plane at slightly different altitudes.
A tight time constraint and possibly large target chains are considered.
The aim is to minimize the overall mission $\Delta V$, while performing a complete tour of the target set. 

Several methods have been proposed to design ADR missions.
Exhaustive, brute-force, approaches \cite{braun2013active}  
and branch and bound search \cite{cerf2013multiple}
have been attempted first.
However, the effectiveness of those methods is limited to small sets of targets, due to the \textit{course of dimensionality}.\cite{lenstra1981complexity} 
Meta-heuristic approaches have thus rapidly gained in popularity, as they allow to find a sub-optimal solution in a reasonable, limited amount of time.
Methods involving constructive heuristics, such as 
Beam Search\cite{izzo2015evolving} and Ant Colony Optimization\cite{shen2017optimization}  have been widely exploited for both ADR 
and asteroids exploration mission planning. An hybridization of the two has also been discussed.\cite{simoes2017multi}
These methods attempt at chaining one target after the other, into a sequence,    
evaluating many (if not all) of the possible ``branches'' departing from a given starting ``node''.
In order to reduce the overall computational effort, 
simplistic transfer models are used for a fast evaluation of any target-to-target transfer cost;
moreover, pre-determined encounter epochs are commonly employed. 
Despite being rather flexible, this formulation usually underperforms in case a complete tour is searched for. 
Instead, meta-heuristic methods,
that rely on an iterative refinement of a set of (possibly random) initial solutions guarantee the tour completeness at any point in the optimization process.
Prominent examples are Tabu Search (TS),\cite{gurtuna2003orbit}
Genetic Algorithm (GA) \cite{federici2019impulsive} 
and Simulated Annealing (SA).\cite{cerf2015multiple}

Among them, Simulated Annealing is one of the most successful, as inherently designed for solving hard combinatorial optimization problems.\cite{Kirkpatrick1983}
Moreover, it has the advantage of being easy to implement and understand, and has solid theoretical foundations, as opposed to other meta-heuristic approaches.
The performance of SA are highly related to the encoding used to map the set of design variables to a configuration of the problem, and consequently to the
neighbor generation functions, or 
\textit{moves}, that create a new solution from the current one.\cite{Bentner2001} 
Therefore, 
a permutation encoding is proposed in this paper for modelling an ADR mission as a TDTSP. \textit{Moves} are randomly chosen within a set of permutation-preserving operators, which guarantee the mutual exchange of cities in the tour without the need for ``repair actions'', whose effectiveness is doubtful.

The paper is organized as follows.
First, formulations for both the classical and the time-dependent TSP are recalled.
Subsequently, a generic ADR mission problem is stated and a bi-level optimization procedure is proposed.
Analogies between the TSP and combinatorial features of ADR missions are highlighted, 
leading to the formulation of a permutation optimization problem.
Simulated Annealing is used for solving the defined discrete optimization problem, which concerns the selection of the ADR mission target sequence and a possibly rough estimates of the encounter epochs.
Once the best sequence and the approximate encounter times have been determined, the whole transfer is optimized through a multi-population self-adaptive Differential Evolution algorithm. Numerical results for sets comprising up to 20 targets are presented.
A conclusion section ends the paper.

\section{The Time-Dependent Traveling Salesman Problem }




\newcommand{\mysum}[3]{\sum_{#1 = #2}^{#3}}

This section briefly reviews two formulations for the 
traveling salesman problem and its time-dependent variants, the former
based on a binary encoding, that leads to an Integer Linear Programming (ILP) problem, the latter on a permutation encoding, that leads to a Permutation Optimization (PO)  problem.

Integer Linear Programming is a rather powerful modeling framework that provides great flexibility for expressing discrete optimization problems. 
Well-established solution algorithms exist that are routinely used for this kind of problems and exploit a combination of linear programming relaxation, branch-and-bound, and branch-and-cut techniques.
Most researches focused on devising ``stronger'' ILP formulations for the analyzed problems, that are,
formulations where the feasible solution set is closer to the feasible set of the corresponding linear programming relaxation.
However, the effectiveness of this formulation is limited by the possibly large amount of memory required by resolving algorithms to deal with medium and large problem instances.

The second approach exploits the more intuitive permutation encoding.
The resulting formulation usually involves a significantly lower number of variables, which allows the use of meta-heuristic methods
capable of reaching good-quality solutions even in larger-problem instances in a limited amount of time.
On the other hand, convergence towards the global optimum cannot be guaranteed,
and often a problem-dependent tuning of the resolving algorithm parameters is required in order to attain satisfying results.

\subsection{Traveling Salesman Problem}

The classical Traveling Salesman Problem deals with the determination of the minimum-distance tour
for a salesman that, departing from his home city, has to visit exactly once each city on a given list, and then come back home \cite{lawler1985traveling}.
Hereafter, let $n$ be the total number of cities, i.e. the cardinality of the cities set $V = \{1, \dots, n\}$.
The cost of going from city $i$ to city $j$ is $c_{ij}$, hence $C = (c_{ij})$ is a $n \times n$ cost matrix.

\subsubsection{ILP Formulation.}

Let $x_{ij}$ be a binary variable indicating if the salesman goes directly from city $i$ to city $j$ $(x_{ij}=1)$ or not ($x_{ij}=0$).
A traditional TSP formulation requires to find the matrix $X = (x_{ij})$ that solves:
\begin{align}
	\min_X \mysum{i}{1}{n} &\mysum{j}{1}{n} c_{ij} x_{ij} &&	&&\label{eq:TSP_cost}\\
    \text{s.t.}&
    \mysum{j}{1}{n} x_{ij} = 1 \text{, }   & i &= 1, \dots, n \label{eq:each_row_con} \\
	\phantom{s.t.} &\mysum{i}{1}{n} x_{ij} = 1 \text{, }   &j &= 1, \dots, n  &&\label{eq:each_col_con} \\
	\phantom{s.t.}	&\sum_{i \in S} \sum_{j \in S} x_{ij} \leq |S| - 1
	\text{, }  & \forall &S \in V,\,|S|\ne 0  &&
	\label{eq:each_subtour}
\end{align}

Equations \eqref{eq:each_row_con} and \eqref{eq:each_col_con} ensures that the salesman departs from (respectively, arrives in) each city exactly once. 
Equation~\eqref{eq:each_subtour} prevents the formation of \textit{subtours}, i.e., closed loops of a subset $S \subset V$ with $|S| < n$ cities.
This requires up to $2^n - 2$ additional constraints, since this condition must hold  for every nonempty subset $S$ of $V$.
As a result, a formulation involving $n^2$ binary decision variables and $\mathcal{O}(2^n)$ 
linear constraints is attained.


\subsubsection{PO Formulation.}

The TSP can be more concisely formulated  as a Permutation Optimization problem. The goal is to find the shortest Hamiltonian cycle of $n$ cities, i.e. the permutation $\bfp=\left\{ p_1, p_2, \ldots, p_{n} \right\} \in \mathcal{P}^n$ of $n$, non repeated, positive, integer elements 
$\{1,\dots,n\}$, which solves: 
\begin{equation}
	\min_{\bfp \in\mathcal{P}^{n}} \mysum{i}{1}{n-1} c_{p_i p_{i+1}} + c_{p_n p_1}
	\label{eq:TSP_cost}
\end{equation}
where $\mathcal{P}^n$ is the set of all permutations of $n$ integers.

\subsection{Time-Dependent TSP (A) }
 
The Time-Dependent TSP (TDTSP) is a generalized version of the classical TSP where the cost of any transfer changes according to time.
Several variants of the problem are reported in literature. Most, if not all, of them can be essentially traced back to two basic cases, namely \textit{A} and \textit{B}.
In the TDTSP-A, the cost of any transfer also depends on its position in the whole tour,\cite{picard1978time} while in the variant TDTSP-B the cost function depends on both departure time and transfer duration. 

\subsubsection{ILP Formulation.}
Several ILP formulations of the TDTSP-A problem has been devised by many researchers in the last decades; a 3-index formulation is here reported.\cite{fox1980n} 
Let us consider $n$ time periods, indexed by $t$, and a time-variant distance tensor $C =( c_{ijt})$, where $c_{ijt}$ indicates the cost of the transfer from city $i$ to $j$ performed in period $t$.
Define $x_{ijt}$ to be a binary variable indicating whether ($x_{ijt} = 1$) or not ($x_{ijt} = 0$) the salesman goes directly from city $i$ to city $j$ in period $t$.
The TDTSP-A is solved by the tensor $X = (x_{ijt})$ which minimizes the tour cost:
\begin{equation}
	\mysum{i}{1}{n}\mysum{j}{1}{n}\mysum{t}{1}{n} c_{ijt} x_{ijt}
	\label{eq:costTDTSPA}
\end{equation}
subject to the following constraints:
\begin{align}
	\mysum{j}{1}{n} \mysum{t}{1}{n} x_{ijt} &= 1 \text{, } \hspace{10mm} i = 1, \dots, n \label{eq:each_row_con_tdtsp} \\
	\mysum{i}{1}{n} \mysum{t}{1}{n} x_{ijt} &= 1 \text{, } \hspace{10mm} j = 1, \dots, n \label{eq:each_col_con_tdtsp} \\
	\mysum{i}{1}{n} \mysum{j}{1}{n} x_{ijt} &= 1 \text{, } \hspace{10mm} t = 1, \dots, n \label{eq:each_time_con_tdtsp} \\
	\mysum{j}{1}{n} \mysum{t}{2}{n} t x_{ijt} - \mysum{j}{1}{n} \mysum{t}{1}{n} t x_{jit} &= 1 \text{, } \hspace{10mm} i = 2, \dots, n \label{eq:flux_con_tdtsp}
\end{align}
Constraints \eqref{eq:each_row_con_tdtsp} and \eqref{eq:each_col_con_tdtsp} ensure that exactly one transfer is performed from and to each city, respectively. Equation~\eqref{eq:each_time_con_tdtsp} establishes that each transfer takes place at a different time period.
Constraint \eqref{eq:flux_con_tdtsp} guarantees that the difference between the departure time and the arrival time referred to a same city is unitary; it also forces the tour to start and end in city $1$.
%

The proposed formulation contains $n^3$ variables and $4n - 1$ constraints. 
An alternative, more compact, $n$-constraints formulation has been also proposed,\cite{fox1980n} where constraints~\ref{eq:each_row_con_tdtsp} to \ref{eq:each_time_con_tdtsp} are replaced by:
 \begin{equation}
	\mysum{i}{1}{n}\mysum{j}{1}{n}\mysum{t}{1}{n} x_{ijt} = n
	\label{eq:total_transfers_tdtsp}
\end{equation}
 
Both formulations can also be adopted to solve the classic TSP,
by assuming a cost tensor $C = (c_{ijt})$  with
$c_{ijt} = c_{ij}\,\, \forall t = 1, \dots, n$, thus
allowing to solve the problem without the need to introduce the subtour elimination constraints.

\subsubsection{PO Formulation.}

The TDTSP-A can be also formulated as a PO problem. 
The objective is to determine the permutation $\bfp \in \mathcal{P}^{n}$ of the numbers $\{1,\dots,n\}$ which solves:
\begin{equation}
	\min_{\bfp \in \mathcal{P}^{n}} {\mysum{i}{1}{n-1} c_{p_i p_{i+1} i} + c_{p_n p_1 n}}
	\label{eq:TDTSPA_cost}
\end{equation}
where the distance tensor $C = (c_{ijt})$ is the same as before.
Any starting city $s$ may be arbitrarily fixed simply by imposing $p_1 = s$.


This formulation does not require any additional decision variable with respect to the one proposed for the classical TSP,
highlighting the simplicity and versatility of
the use of a permutation encoding. 

\subsection{Time-Dependent TSP (B) } 

A more general time-dependent TSP regards a scenario where 
a given amount of time is given to complete the tour, and 
the cost for moving from a city to another depends on
both  starting time and duration of the transfer itself.
A finite number $n_t \ge n+1$ of possible departure/arrival instances is assigned, usually far larger than the number of cities to visit, giving the salesman some additional degrees of freedom when designing his trip. 


\subsubsection{ILP Formulation.}

Let us consider $n_t$ times, indexed by $t$, and a distance tensor $C = (c_{ijtm})$, where $c_{ijtm}$ is the cost of departing from city $i$ at time $t$ and arrive at city $j$ at time $t + m$. 
Again, binary decision variables $x_{ijtm}$ are correspondingly defined.
Formally, the problem can be posed as:
\begin{align}
    \min_X &\mysum{i}{1}{n}\mysum{j}{1}{n}\mysum{t}{1}{n_t}\mysum{m}{1}{M} c_{ijtm} x_{ijtm}\\
    \text{s.t.}		&\mysum{i}{1}{n}\mysum{j}{1}{n}\mysum{t}{1}{n_t}\mysum{m}{1}{M} x_{ijtm} = n \label{eq:total_transfers_m_tdtsp} \\
    			&\mysum{j}{1}{n}\mysum{t}{1}{n_t}\mysum{m}{1}{M} x_{ijtm}  = 1 \text{, } & 
    			i = 1, \dots, n \label{eq:each_row_con_m_tdtsp} \\
	&\mysum{i}{1}{n}\mysum{t}{1}{n_t}\mysum{m}{1}{M} x_{ijtm}   = 1 \text{, } & 
	j = 1, \dots, n \label{eq:each_col_con_m_tdtsp}\\
	&\mysum{j}{1}{n}\mysum{t}{2}{n_t}\mysum{m}{1}{M} t x_{ijtm} - \mysum{j}{1}{n}\mysum{t}{1}{n_t}\mysum{m}{1}{M} t x_{jitm} = \mysum{j}{1}{n}\mysum{t}{1}{n_t}\mysum{m}{1}{M} m x_{jitm} \text{, } & 
	i = 2, \dots, n \label{eq:flux_con_m_tdtsp}\\
	&\mysum{i}{1}{n}\mysum{j}{1}{n}\mysum{t}{1}{n_t}\mysum{m}{1}{M} m x_{ijtm} \leq n_t \label{eq:max_m_sum}
\end{align}
where $M \le n_t - (n - 1)$ is an \textit{a priori} assigned value for the maximum duration of all possible transfers.
Constraint~\eqref{eq:total_transfers_m_tdtsp} guarantees a total number of transfers equal to $n$,
%
constraints~\eqref{eq:each_row_con_m_tdtsp} and \eqref{eq:each_col_con_m_tdtsp} ensure one departure and one arrival for each city, constraint~\eqref{eq:flux_con_m_tdtsp} enforces that every transfer has length $m$ and constraint~\eqref{eq:max_m_sum} limits the total duration of the tour.
%

\subsubsection{PO Formulation.}

An ingenious formulation allows to pose the TDTSP-B as a PO problem,  
by using as decision variables a permutation of the available departure times from the cities.
%
%
Let $\bftau=\left\{1, \ldots, n_t \right\}$ be a discrete time grid with $n_t$ departure/arrival instances.
Both the city sequence and the location of the departure epochs over the grid are needed for evaluate the tour cost. However, all decision variables can be collapsed into ``one'' decision variable, 
that is a permutation $\bfPi \in\mathcal{P}^{n_t-1}$, with a number of elements equals to the (discrete) number of available departure epochs minus one.
Permutation elements with a value greater than $n$ are considered as blanks, thus revealing the city sequence $\bfp$,
while the position of non-blank elements reveals a vector of departure epochs $\bft=\left[t_1, t_2, \ldots, t_{n} \right]$,
where the element $t_k$ is the departure epoch from city $p_k$.
The use of a permutation with $n_t-1$ elements instead of exactly $n_t$ allows to guarantee that the final arc which closes the tour last at least on time unit.



The TDTSP-B can be thus written with a PO formulation as:
\begin{align}
\label{eq:TimeDiscreteProblem}
\min_{\bm{\Pi} \in \mathcal{P}^{n_t-1}} & \sum_{k=1}^{n-1}{c_{p_{k} p_{k+1} t_{k} (t_{k+1} - t_{k})} + \bar c_{p_{n} p_{1} t_{n}}}  \\
\bfp =& \left\{\Pi_h \mid \Pi_h\leq n,\,\, \forall h \in [1, {n_t - 1}] \right\}\\
\bft =& \left\{h  \mid  \Pi_h\leq n,\,\, \forall h \in [1, {n_t - 1}]\right\} 
\end{align}

\noindent
where the distance tensor $C = (c_{ijtm})$ is the same as for the ILP formulation, 
and 
$$\bar c_{ijt} = \min_{r\in[t+1,\, n_t]}{c_{ijt(r - t)}}$$
is the minimum viable return cost for closing the tour, given the remaining time intervall.




\section{ADR mission statement}\label{sez:statement}

Let us consider a 
set of $N$ prescribed target debris that move on circular coplanar orbits at slightly different altitudes under a keplerian dynamical model. 
For each target body $A_i$,  orbital radius $r_{A_i}$ and right ascension at starting time $\theta_{A_i}(t_0)$ are assigned. 
The velocity is constant and equals to $v_{A_i}=\sqrt{{\mu}/{r_{A_i}}}$, while the angular position at any time is given by $\theta_{A_i}(t) = \theta_{A_i}(t_0) + \sqrt{{\mu}/{r_{A_i}^3}}\,(t-t_0)$,
where $\mu$ is the gravitational parameter of the central body.
The chaser is also assumed to be initially on a circular orbit of radius $r_0$ at the right ascension $\theta_{0}(t_0) = 0$, on the same orbital plane as all the targets. The problem is thus planar.
The goal is to design an impulsive multi-rendezvous transfer trajectory that allows the chaser to perform a complete tour of the prescribed set of targets within a specified maximum time-length of the entire mission $T_M$, minimizing the overall mission $\Delta V$.
 

%

\begin{figure} [!h]
    \centering
    \includegraphics[width = 0.55\textwidth]{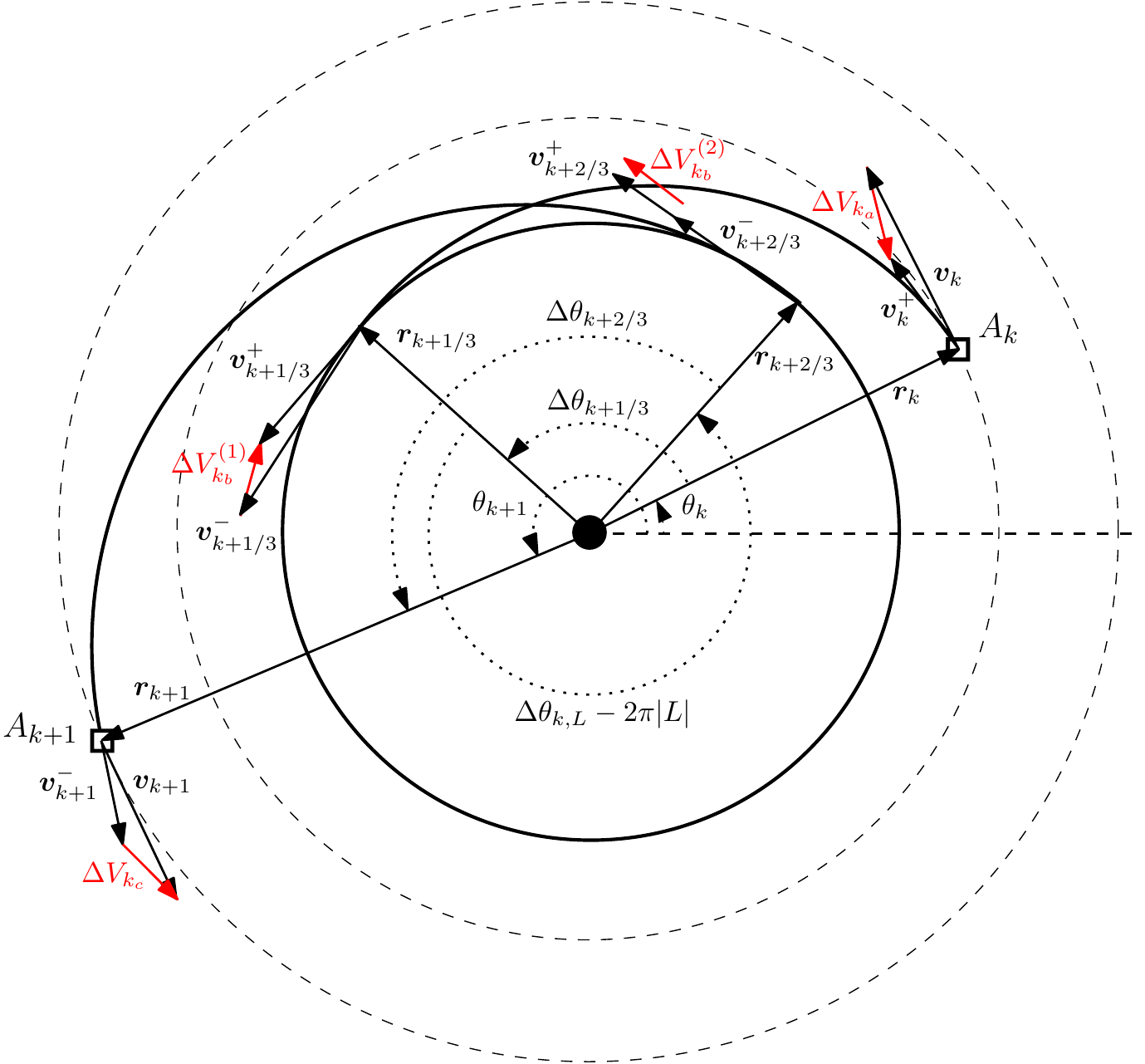}
    \caption{Trajectory sketch for the $k$-th leg.}
    \label{fig:Leg}
\end{figure}

Let us assume that a sequence $S_{A} = \left\{A_1,\,A_2,\,\ldots,\,A_N \right\}$ of N non-repeating bodies to encounter and the corresponding set of (monotonically increasing) encounter times $\bft = \left\{t_1,\,t_2,\,\ldots,\,t_N \right\}$, with $t_N = T_M$,
have been assigned,
so that the integer $A_j \in [1, N]$ identifies the target met at time $t_j$. 
%
%
%
%
The overall trajectory of the chaser can be decomposed into a series of target-to-target body legs. 
The $k$-th leg departs from body $A_{k}$ (with $A_0 = 0$ denoting the chaser) at time $t_{k}$ and arrives at the body $A_{k+1}$ at time $t_{k+1}$, for $k = 0,\dots,N-1$.  


The rendezvous condition requires that, at the ending point of the leg, position and velocity of the chaser are the same as the target:
\begin{align}
    \bfr(t_k) &= \bfr_{A_k}(t_k) = \bfr_k &\qquad \forall\, k \in \left\{0,1,\dots,N\right\}\\
    \bfv(t_k) &= \bfv_{A_k}(t_k) = \bfv_k &\qquad \forall\, k \in \left\{0,1,\dots,N\right\}
\end{align}

Being four
the maximum number of impulses for an optimal time-constraint planar rendezvous,\cite{prussing1967optimal} each body-to-body transfer leg is made up of a sequence of three ballistic \textit{arcs}, named ``\textit{a}'', ``\textit{b}'', and ``\textit{c}'', joined by impulsive maneuvers located  at the departure, at the two internal points labeled with subscripts ``$k+1/3$'' and ``$k+2/3$'', and at the arrival point, respectively.

A position formulation is here considered, that is, the trajectory is parameterized with respect to radii $r_{k+1/3}$, $r_{k+2/3}$ and anomalies $\theta_{k+1/3}$, $\theta_{k+2/3}$ at the internal maneuvering points. 
Spacecraft velocities immediately before $\bfv_{k+1/3}^-$, $\bfv_{k+2/3}^-$ or after $\bfv_{k+1/3}^+$, $\bfv_{k+2/3}^+$ the maneuvers are found by solving either a geometrical problem or a Lambert problem.

\subsection{Geometrical Problem}
Let us  first consider an arc ``\textit{a}'' connecting the points ``$k$'' and ``$k+1/3$''.
Two families of ellipses that connects  $\bfr_k$ and $\bfr_{k+1/3}$ exist.
They might be parameterized as a function of the semi-major axis and  
labeled as \textit{fast} and  \textit{slow} families.\cite{cacciatore2008optimization}

\begin{figure}[!h]
 	\centering
 	\subfigure[Transfer ellipse flight time $\Delta t$ and semi-major axis $a$\newline as a function of the $y$ parameter.]
 	{\label{fig:MonteA}
 		\includegraphics[width=0.45\textwidth]{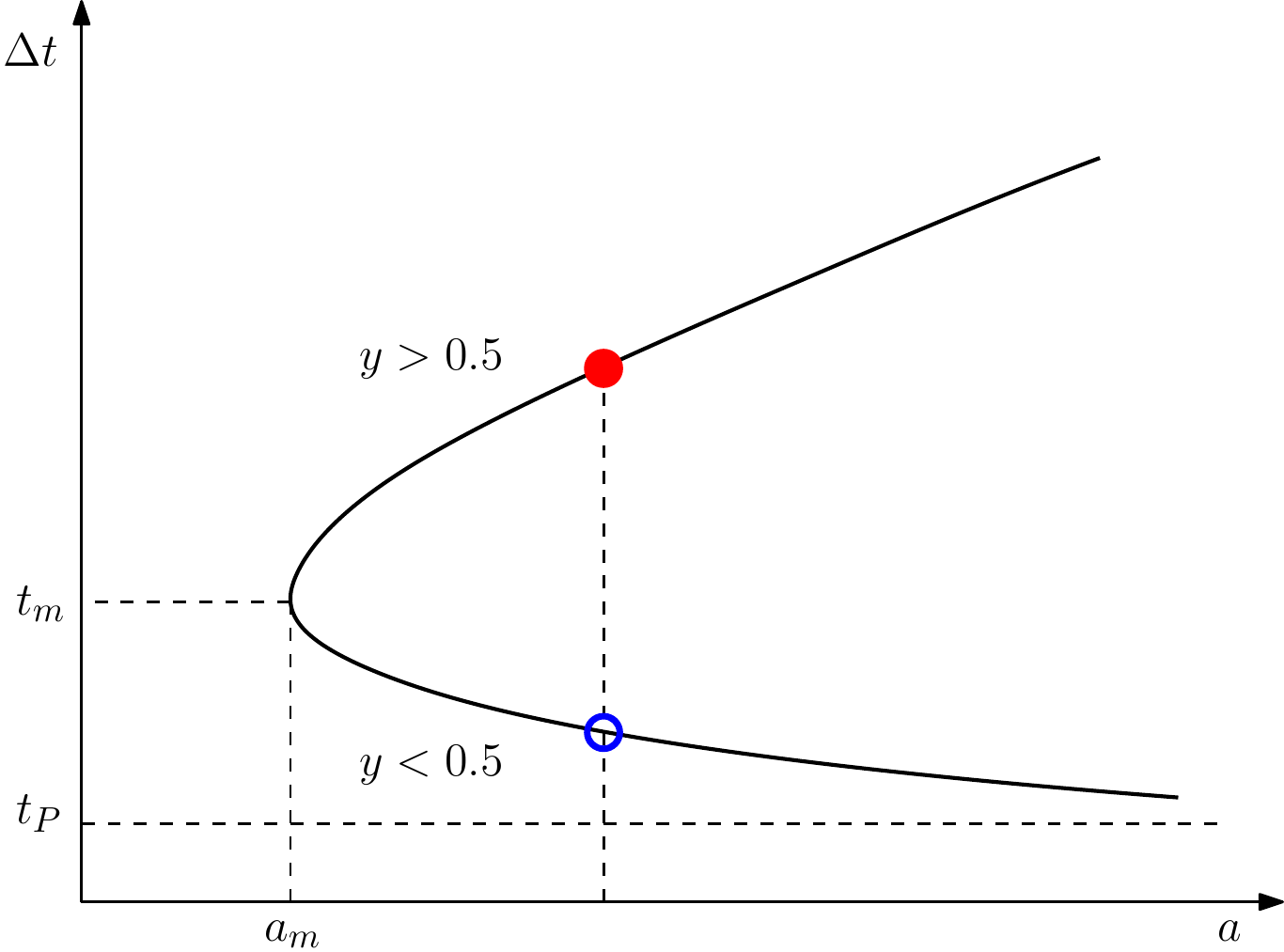}} 
 	\subfigure[ Geometrical construction of the transfer ellipse.]
 	{\label{fig:yarc}
 		\includegraphics[width=0.375\textwidth]{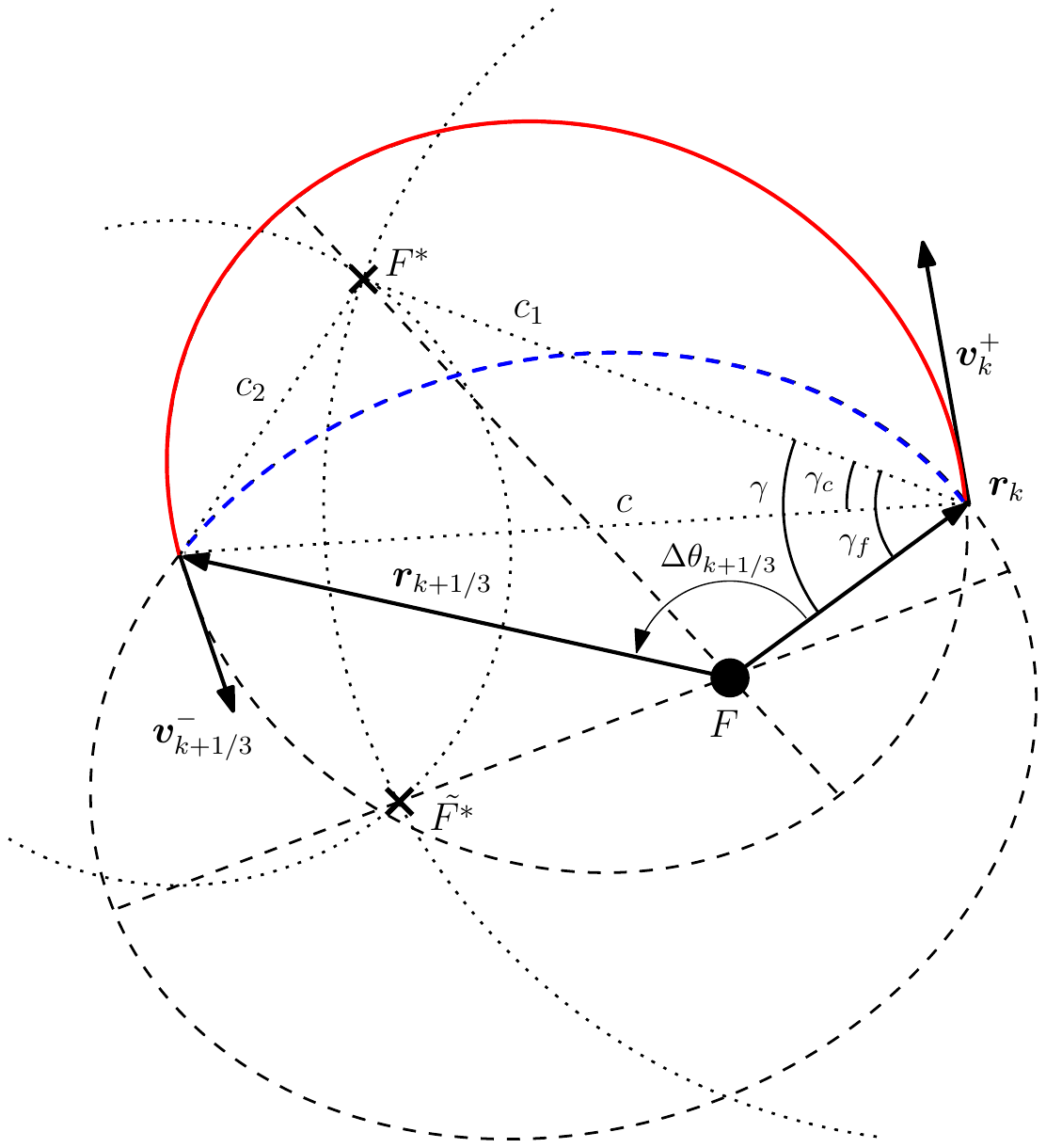}}
 	\caption{Geometrical Problem.  
 	Filled marker refers to the slow solution, corresponding to the solid-line transfer arc; empty marker refers to the fast, dashed-line solution.
 	}
 	\label{fig:yProblem}
\end{figure}

\noindent 
Let us introduce a non-dimensional parameter $y\in[0,1]$ so that:
\begin{equation}
    \label{eq:yk}
    a = \frac{a_m}{4y\,(1 - y)}
\end{equation}
where $a_m = (r_{k} + r_{k+1/3} + c)/4$ is the smallest semi-major axis which connects the ending points, and $c =\norm{\bfr_{k}-\bfr_{k+1/3}}$ is the chord distance between them.
%
Pairing the fast solutions with values $y<1/2$,
and the slow solutions with values $y>1/2$, the elliptic arc connecting the two assigned endpoints is uniquely identified for any given choice of $y$. 
Figure~\ref{fig:yarc} depicts the 
geometric construction that allows an unambiguous definition of the transfer ellipse, that is, semi-major axis $a$, eccentricity $e$  and argument of pericenter $\omega$,
as soon as $y$ and $\Delta \theta$ are known.
Chord lengths $c_1=2 a-r_k$, and $c_2 = 2 a - r_{k+1/3}$ are evaluated first.
Angle $\gamma$  follows from:
\begin{equation}
\gamma  = \acos{\left( \frac{r_{k}^2 - r_{k+1/3}^2 + c^2}{2 r_{k}\, c} \right)}    
\end{equation}

hence $\gamma_f = \gamma - \gamma_c$, where:
\begin{align}
\gamma_c &=
 \begin{cases}
    \acos{\left( \dfrac{c_1^2 - c_2^2 + c^2}{2 c_1 c} \right)} &\mbox{if}\,\,\, y \leq 0.5\\
    - \acos{\left( \dfrac{c_1^2 - c_2^2 + c^2}{2 c_1 c} \right)} &\mbox{if}\,\,\, y > 0.5
    \end{cases} 
\end{align}    
Eventually, the eccentricity is found as
\begin{equation}
    e = \frac{\sqrt{c_1^2 + r_{k}^2 - 2 c_1\, r_{k} \cos{\gamma_f}}}{2 a}
\end{equation}
Velocities at both endpoints ($\bfv_{k}^+$ and $\bfv_{k+1/3}^-$) follow from standard equations of the two-body problem. 
Transfer time $\Delta t_{a}$ can also be evaluated and, 
consequently, the epoch at the intermediate maneuver $t_{k+1/3} =t_k+ \Delta t_{a}$ is obtained.
The same procedure holds for the arc connecting points ``$k+2/3$'' and ``$k+1$'', resulting in a similar geometrical definition of the third arc ``\textit{c}''. As a result, one has
\begin{equation}
 \left[\bfv_{k}^+,\, \bfv_{k+1/3}^- \right]\leftarrow yArc\left(\bfr_k, \bfr_{k+1/3},y_{k,a}\right)   
\end{equation}
\begin{equation}  
 \left[\bfv_{k+2/3}^+,\, \bfv_{k+1}^- \right]\leftarrow yArc\left(\bfr_{k+2/3}, \bfr_{k+1},y_{k,c}\right) 
\end{equation}

\noindent
The cost of the maneuvers at the departure and arrival points are thus evaluated as:
\begin{align}
    \Delta V_{k_a} &= \norm{\bfv_{k}^+ - \bfv_{k}}\\
    \Delta V_{k_c}& = \norm{\bfv_{k+1} - \bfv_{k+1}^-}
\end{align}

\subsection{Multi-revolution Lambert problem}
The central arc ``\textit{b}'', which connects points ``$k+1/3$'' and ``$k+2/3$'', cannot be dealt with in the same fashion. In fact, for a given choice of the parameters $y_{k,a}$ and $y_{k,c}$ the maneuvering epochs $t_{k+1/3}$ and $t_{k+2/3}$, and, consequently, the travel time, are assigned. 
A multi-revolution Lambert problem 
can be formulated, being  the position vectors $\bm{r}_{k+1/3}$, $\bm{r}_{k+2/3}$ and the travel time $\Delta t = t_{k+2/3} - t_{k+1/3}$ known.
This problem admits $1+2 n_{max}$ solutions, where $n_{max}$ is the maximum allowed number of revolutions:
one solution for the 0-revolution transfer arc and two additional solutions, namely left and right branch, for each $n$-revolution transfer orbit.
Let us introduce an integer parameter  $L\in[-n_{max},n_{max}]$
indicating the solution corresponding to the $\left|L\right|$-revolution transfer orbit
and side $sign(L)$, positive for the right branch, negative for the left one.
For the $L$-th solution, 
the velocity vectors immediately after the second impulse $\bm{v}_{k+1/3}^+$ and just before the third one $\bfv_{k+2/3}^-$ can be evaluated as
according to the algorithm by Izzo~\cite{izzo2015revisiting}, that is:
\begin{equation}
   \left[\bm{v}_{k+1/3}^+, \bm{v}_{k+2/3}^- \right]\leftarrow Lambert\left(\bm{r}_{k+1/3}, \bm{r}_{k+2/3}, t_{k+2/3} - t_{k+1/3}; L \right)
\end{equation}
Hence, the total cost of the two internal maneuvers is:  
\begin{align}
    \Delta V_{k_b} &= \norm{\bm{v}_{k+1/3}^+ - \bm{v}_{k+1/3}^-} +
      \norm{\bm{v}_{k+2/3}^+ - \bm{v}_{k+2/3}^-}
\end{align}
 
Instead of treating $L$ as an optimization variable, an enumeration approach could be used, that is, all the $2 n_{max}+1$ possible solutions are computed and the one with the lowest total $\Delta V$ is chosen.
However, as we are considering transfers between close orbits, we can make an educated guess and safely restrict the analysis to just three scenarios, that is the chaser performs the same number of revolution as the target, one more, or one less,
corresponding to six solutions (three right and three left branches).

\subsection{MINLP formulation}
According to the proposed formulation, the overall trajectory can be parameterized by using a set of $8\, N$ parameters, that is:
\begin{equation}
    \bm{x} = \bigcup\limits_{k=1}^N \bm{x}_k
\end{equation}
where:
\begin{equation}
\bm{x}_k = 
\left\{ A_k,\, t_k  \right\}
\cup
\left\{ r_{k+1/3},\,\Delta \theta_{k+1/3},\, y_{k,a}  \right\}
\cup
\left\{ r_{k+2/3},\,\Delta \theta_{k+2/3},\, y_{k,c}  \right\}\,
\label{eq:x}
\end{equation} 
%
%
Eventually, the impulsive time-constrained MRR optimization problem can be 
formulated as:

\begin{equation}
\label{eq:MINLP}
\begin{array}{ll}
P & = 
\begin{cases}
\min\limits_{\bm{x}} \Delta V_{tot}(\bm{x}) \\
\mbox{s.t.}\quad \bm{x}_L \leq \bm{x} \leq \bm{x}_U\\
\end{cases}
\end{array}
\end{equation}
where the overall cost of the MRR trajectory is:
\begin{equation}
\Delta V_{tot} = \sum_{k = 0}^{N-1}{\left(\Delta V_{k_a} +\Delta V_{k_b} +\Delta V_{k_c} \right)
}
\end{equation}
and $\bm{x}_L, \bm{x}_U$ are the lower and upper bounds of the design variables, respectively.
 This problem
involves the simultaneous optimization of both integer variables (defining the encounter sequence) and real-value decision variables (such as, radius and anomaly at the maneuvers) and it is thus labeled as a  Mixed-Integer Nonlinear Programming (MINLP) problem.

\section{Bi-Level Optimization Approach}\label{sez:bilevelOpt}

The MINLP problem in Eq.~\eqref{eq:MINLP} belongs to
the class of NP-hard problems, hence
no  deterministic  algorithm exists for finding the optimal solution in polynomial-time. 
A variety of stochastic meta-heuristic techniques have been developed over the last decades aiming at 
attaining a (often sub-optimal) good-quality solution in a reasonable, limited amount of time.
However, as the problem dimension increases, the required computational time may become prohibitive.
%

Instead of solving the problem \textit{as a whole}, 
one  might attempt to decompose the problem into simpler sub-problems, that could be (more or less easily) solved separately, and eventually their solutions can be recomposed into the original problem solution.
%
For the problem at hand, a bi-level approach can be pursued, by isolating
i) an outer level that concerns the definition of the encounter sequence and a (possibly rough) evaluation of the epochs at each encounter, while details of each  body-to-body transfer leg are neglected; 
ii) an inner level which deals with the optimization of each body-to-body transfer with full details, 
assuming that departure and arrival bodies are assigned; encounter epochs may or may not be fixed.
%

The two layers are interconnected: 
the outer layer requires a measure of the cost associated to each transfer leg for ``weighting'' the quality of a certain encounter sequence,
even though the actual $\Delta V$ of each leg can be evaluated only by solving the full-transfer optimization problem, that is, the inner-layer problem. 
On the other hand, the inner layer requires the definition of the encounter sequence and rendezvous epochs, which, in turn, are the output of the combinatorial, outer-layer problem.
In practice,  the two problems might be solved sequentially provided that a way, that is, a heuristic, 
exists for attaining a reasonable estimate of the transfer cost without solving the full optimization problem.
Once the heuristic has been established,
the outer-level combinatorial problem, or touring problem, is isolated and solved first; its solution is then used as initial guess for the inner-level problem.
%
%
%

\subsection{Cost Estimate for a Single Rendezvous Leg} \label{sez:AnaliticalSolution}

This section presents an analytical, sub-optimal, four-impulse strategy to assess the $\Delta V$ of a trajectory leg, for any assigned pair of 
departure and arrival bodies that fly on circular orbits,
which fairly approximates the behavior of the time-fixed optimal solution,\cite{Zavoli2015,ZavoliAlaska}
when the allowed travel time is sufficiently large.


\begin{figure}[h]
	\centering
	\includegraphics[width=0.4\textwidth]{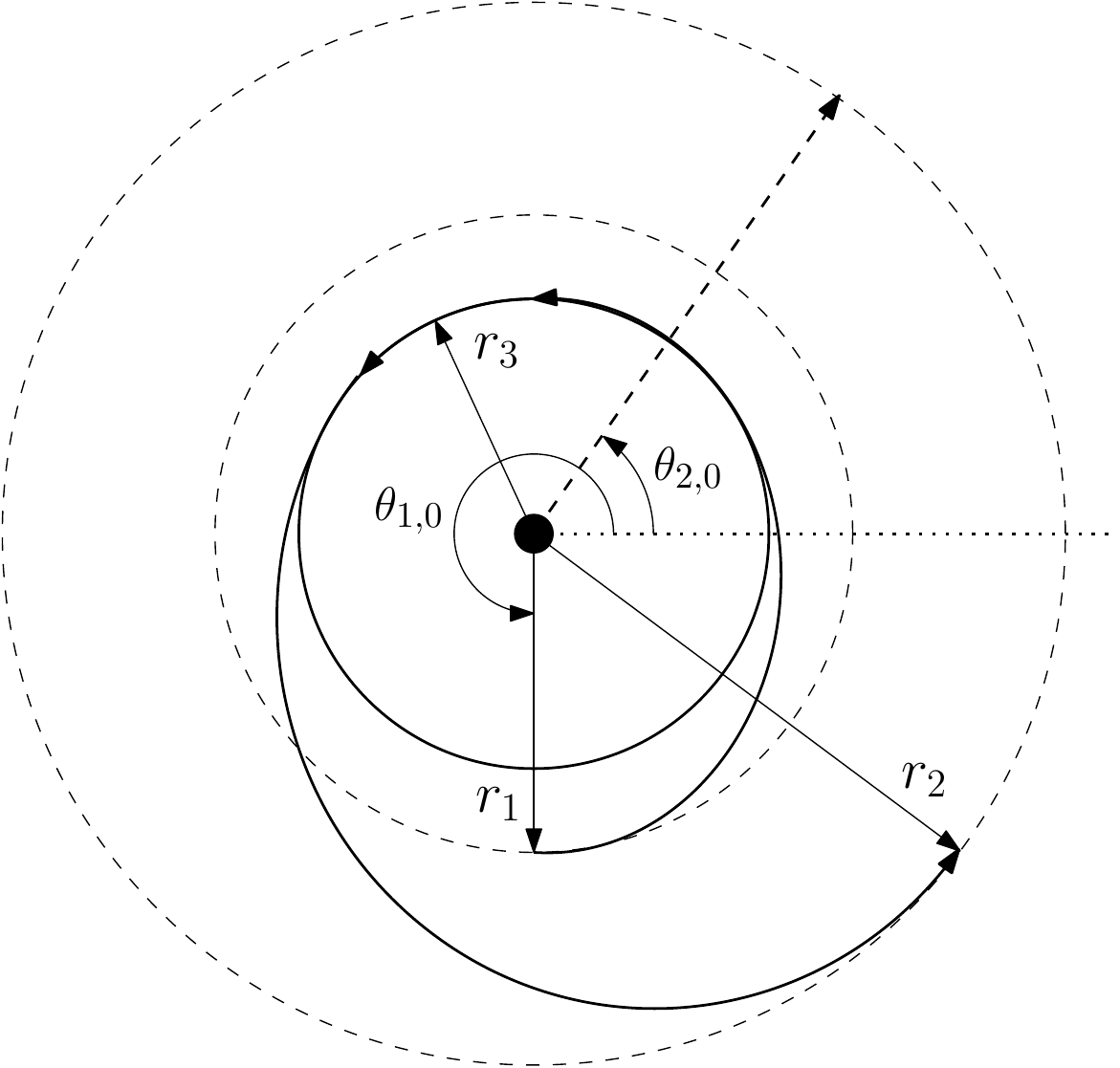}%
	\caption{Chaser trajectory according to the sub-optimal rendezvous strategy adopted as heuristic.}
	\label{fig:trajAnalitic}
\end{figure}

Assuming that departure and arrival orbits are not too far apart,
the minimum-$\Delta V$  solution is represented by a Hohmann transfer,
possibly preceded and/or followed by coasting arcs on the departure and/or arrival orbit 
which allow for the correct phasing required by this kind of maneuver.
Let $\theta_{1,0}$ (respectively, $\theta_{2,0}$), be the true anomaly at time $t=0$ of the departure (respectively, arrival) body, flying on circular orbits of radius $r_1$ (respectively, $r_2$).
%
%
The departure coasting arc duration, that is, the time $T_{wait}$ required to attain the correct phase $\gamma^\star = \pi-\omega_2 T_{H_{12}}$ between departure and arrival body,  
can be evaluated as: 
%
\begin{equation}
\begin{array}{ll}
     T_{wait} =&
     \begin{cases}
     \dfrac{\theta_{2,0}-\theta_{1,0}-\pi+\omega_2 T_{H_{12}} + 2(z+1)\pi}{\omega_1-\omega_2} & \mbox{if}\quad r_1 < r_2 \\[10pt]
     \dfrac{\theta_{2,0}-\theta_{1,0}-\pi+\omega_2 T_{H_{12}} +2z\pi}{\omega_1-\omega_2} & \mbox{if}\quad r_2 < r_1
     \end{cases}
\end{array}
\label{eq:Twait}
\end{equation}
with $z = \left\lfloor{\dfrac{\pi - \omega_2 T_{H_{12}} + \theta_{1,0} - \theta_{2,0}}{2\pi}}\right\rfloor$.

By comparing the available maximum transfer time $T_{max}$ with the sum of the waiting time $T_{wait}$ plus the time spent on the Hohmann transfer $T_{H_{12}}=\pi\sqrt{a_{12}^3/\mu}$,
one obtains a condition for the availability of the Hohmann transfer:
\begin{equation}
\label{eq:H}
T_{wait}+T_{H_{12}} \le T_{max} 
\end{equation}

Whenever Eq.~(\ref{eq:H}) holds, the  cost of the transfer leg is easily evaluated as $\Delta V_h$ of the Hohmann transfer.
%
If the Hohmann transfer is not possible,
the  mission scheme depicted in 
Figure~\ref{fig:trajAnalitic} is adopted:
the maneuvering spacecraft is injected into a circular (either internal or external) waiting orbit of radius $r_3$ with an Hohmann transfer ``1-3'', 
of semi-major axis $a_{13}=(r_1+r_3)/2$ and
duration $T_{H_{13}}=\pi\sqrt{a_{13}^3/\mu}$, in order to adjust its phase with respect to the target body. A second Hohmann transfer ``3-2'', of semi-major axis $a_{23}=(r_2+r_3)/2$ and
duration $T_{H_{23}}=\pi\sqrt{a_{23}^3/\mu}$ is then used to close the rendezvous.
The rendezvous equation, which imposes the equality of chaser and target position at the end of the maneuver, is enforced:
%
%
\begin{align}
\label{eq:rendezvous}
\Delta\theta_0 =  \left(T_{max}-T_{H_{13}}-T_{H_{23}}\right)\omega_3   -  T_{max} \omega_2  + 2 (1 - k_{rev})\pi
\end{align}

\noindent
with $\Delta\theta_0 = \theta_{2,0}-\theta_{1,0} \in [0,2\pi]$.


This nonlinear equation in $r_3$ admits a family of solutions, parameterized by $k_{rev}\in \mathbb{Z}$,
that represents the additional number of revolutions performed by the maneuvering spacecraft with respect to the target.
%
However, being interested to the minimum $\Delta V$ solution only, one can safely restrict the search to the cases 
corresponding to the largest inner waiting orbit ($k_{rev} = 0$ or $k_{rev} = 1$) and the smallest outer waiting orbit ($k_{rev} = 0$ or $k_{rev} = -1$).
Depending on the initial relative phasing of the two bodies, the relative angular velocity, and the maximum allowed travel time, each of these 
solutions could be the best one. 
So, the three solutions are evaluated and compared each other. The one with the lower cost is retained and the corresponding $\Delta V$ is used as a cost estimate.

\subsection{Outer-Level Optimization}\label{sez:Tour}

This section presents three formulations of the 
outer-level ADR mission design problem, which match the permutation-based formulations of the TSP, TDTSP-A and TDTSP-B previously discussed.

In order to carry on the analogy between an ADR mission and the TSP, let us assume that: (i) in an ADR mission, the TSP cities correspond to target debris; (ii) each city-to-city distance is substituted by the velocity variation along the corresponding orbital transfer; so the distance tensor $C$ is replaced by a $\Delta V$ tensor; (iii) the initial city, corresponding to the initial position of the chaser, is preassigned; by setting null the cost of returning to such initial condition, a ``closed'' TSP problem can be still solved; in fact, since the cost of the transfer from the ``final'' city to the pre-assigned initial one is null, it does not affect the objective function.

\subsubsection{Time-Free Tour.}\label{sez:timeFreeTour}

Under the assumption that both departure and arrival bodies fly on close coplanar circular orbits and that the allowed mission time is sufficiently large, 
the optimal transfer is always a Hohmann transfer.
The cost $\Delta V(i,j)$  for moving from a departure body $i$ to an arrival body $j$ is the same, regardless of the specific departure/arrival epochs.
The problem thus reduces to the search for the sequence of encountered bodies that minimizes the total velocity increment: it is clear the parallelism of such problem with the classical TSP.
So as to speed-up the resolving algorithm, transfer costs for any pair of arrival/departure bodies can be preliminary evaluated and collected in a cost matrix $\Delta V_{ij}$ of dimensions $(N+1) \times N$, being $N$ the number of targets.

\subsubsection{Time-Fixed Time-Uniform Tour.}\label{sez:timeUniformTour}
The solution of the previous problem provides a lower bound on the overall tour cost, but may not provide a reasonable guess to the solution of the original problem. In fact, the perfect phasing required by the Hohmann transfer might never occur 
due to the existing time-constraint on the overall mission duration.
For time-fixed rendezvous maneuvers, transfer costs are highly sensitive to the initial phasing, that is, to the  departure epoch, and to the allowed transfer duration.
However, under the assumption that the duration of all transfers is exactly the same, 
a relaxed problem
can be formulated.
The epoch of the $k$-th encounter is readily available as:  $t_k = k\, T_{M}/N$. 
The problem reduces again to the search for the optimal permutation of the target sequence $\bfp\in\mathcal{P}^{N}$,
but, in this case, 
the cost associated to the transfer from a target $i$ to a target $j$ also depends on the position $k$ of the leg into the sequence.
Thus, the \textit{time-fixed time-uniform tour} problem is the same as the TDTSP-A previously described.
Again, all transfer costs can be preliminary evaluated and collected in a 3-dimensional cost tensor $\Delta V_{ijk}$.

\subsubsection{Time-Fixed Time-Discrete Tour.}\label{sez:timeDiscreteTour}


Eventually, if 
the encounter epochs are discretized over a finer time-grid, one obtains a combinatorial optimization problem
with a solution  closer to the solution of the original problem and 
a formulation analogous to that of the TDTSP-B.
However some remarkable differences arise in the formulation, due to the fact that the starting city (chaser) is unambiguously defined, and no return transfer to the chaser orbit is needed. 

\begin{figure}[h!]
	\centering
	\includegraphics[width=0.5\columnwidth]{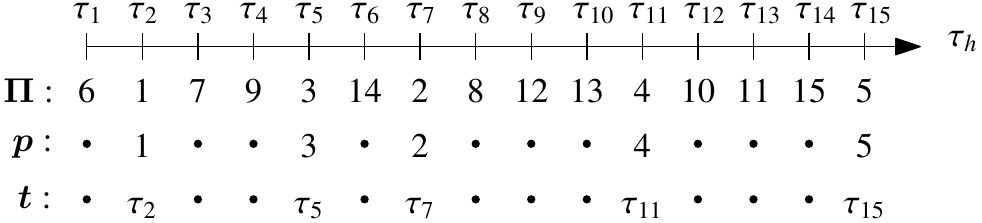}
	\caption{An example of a permutation encoding/decoding for a 5x3 tour.}
	\label{fig:permTD}
\end{figure}

For the sake of clearness, the complete time-fixed time-discrete formulation is here summarized.
Let $\bftau=\left\{\tau_0, \tau_1, \ldots, \tau_{N^s} \right\}$ be a discrete time grid with uniform, equidistant points, 
where $N^s = N D$, $D$ is the number of divisions introduced into each ``previously considered'' time-slot, and $\Delta V(i,j,\tau_h, \tau_{h+m})$ 
denotes the cost for moving from body  $i$ to body $j$, departure epoch $\tau_h$, and arrival epoch $\tau_{h+m}$.
A permutation $\bfPi \in\mathcal{P}^{N^s}$ with a number of elements equals to the (discrete) number of available encounter epochs is used as optimization variable.
In the same fashion as in the TDTSP-B, permutation elements of value greater than $N$ are ignored, 
thus revealing both an encounter sequence $\bfp=\left[p_1, p_2, \ldots, p_{N} \right]$, and an encounter epochs' vector $\bft=\left[t_1, t_2, \ldots, t_{N} \right]$.
The \textit{time-fixed time-discrete tour} problem can now be written as:
\begin{align}
\label{eq:TimeDiscreteProblem}
\min_{\bfPi \in\mathcal{P}^{N^s}} & {\sum_{k=1}^{N}{\Delta V(p_{k-1},p_k,t_{k-1},t_k)}} \\
\bfp =& \left\{\Pi_h \mid \Pi_h\leq N,\,\, \forall h \in [1, {N^s}] \right\}\\
\bft =& \left\{\tau_h = h\Delta T \mid  \Pi_h\leq N,\,\, \forall h \in [1, {N^s}]\right\} 
\end{align}
with $p_0=0$ being the chaser, $t_0 = 0$ the initial time and $\Delta T = T_{M}/(N D)$ the time unit of the time grid.

An encoding/decoding example is proposed in Figure~\ref{fig:permTD} for $N = 5$ and $D=3$, showing a permutation $\bfPi \in \mathcal{P}^{N^s}$ which reveals a target sequence 
$\bfp = \left\{ 1,3,2,4,5\right\}$ and encounter epochs \mbox{$\bft = \left\{\tau_2, \tau_5, \tau_7, \tau_{11}, \tau_{15}  \right\}$}.


It is worthwhile to remark that the number of divisions $D$ should be kept small, because
	i) a rough evaluation of the encounter epochs is sufficient, as
	the attained solution will be further refined within the inner-level optimization step; 
	ii)  a large number of divisions makes the problem too similar to the original MINLP, hence more difficult to solve;
	iii) the proposed heuristic works well if there is enough time to perform several revolutions;  reducing the minimum valid travel time $\Delta T = T_M/(N D)$ makes the heuristic less reliable and may undermine our efforts.
As a result, a number of divisions $D =2$ or $D =3$ appears as a good trade-off value for the problem at hand.
%


As in the other touring problems,
one may pre-compute all transfer costs for speeding up the function evaluation (hence the whole optimization process).
A 4D tensor 
of dimensions 
$\left[ N+1,\, N,\, N D - 1,\,  N(D-1) - 1  \right]$ 
is needed in principle.
However, one may notice the monotonic, non-increasing, behavior of $\Delta V$ with transfer time 
and decide to 
limit the calculus to  transfers of duration $M \Delta t$, 
assuming the same cost for longer transfer. 
A  4D tensor of dimensions  
$\left[ N+1 , N , N D -1 , M \right]$  
would now be required.
Apart from reducing the size of the tensor, this treatment has the  
additional benefit of guiding  the solver toward a trajectory with more uniform travel times,  
which might be good from an operational point of view. 

\subsection{Inner-Level Optimization}\label{sez:Inner}
Assuming that the target sequence  $S_A$ has been selected and a rough estimation of the encounter epochs is available, 
the MINLP problem in Eq.~\eqref{eq:MINLP} reduces to a NLP problem.
%
Two scenarios can be investigated: i) encounter epochs $t_k$ are kept fixed at their nominal values $\bar{t}_k$; ii) encounter epochs are free to be optimized.
In the first case, each body-to-body transfer can be solved independently from the others, thus reducing the (now) $6 N$ problem to solving $N$ easier sub-problems each one of dimension 6,
being $N$ the prescribed number of targets to encounter.
In the second case, encounter epochs may vary in a neighborhood of the reference value, leading to an improved solution, but the whole trajectory must be fully optimized. Lower and upper bounds of the encounter epochs 
are selected so that:
\begin{equation}
t_k \in \left [ \bar{t}_k - \frac{\bar{t}_{k}-\bar{t}_{k-1}}{2}, \bar{t}_k + \frac{\bar{t}_{k+1}-\bar{t}_k}{2} \right]
\end{equation}

\section{Differential Evolution} \label{sez:DE}

In the present paper,
the inner level optimization is carried out by using an 
in-house optimization code named \textit{EOS} (Evolutionary Optimization at Sapienza),
developed in the contest of the Global Trajectory Optimization Competitions\cite{Colasurdo2014190,Casalino2016} 
and previously  applied with success to both unconstrained\cite{federici2018preliminary} and constrained\cite{Federici2019Integrated} space trajectory optimization problems.
\textit{EOS} implements a multi-population self-adaptive Differential Evolution (DE) algorithm, with a synchronous island-model for parallel computation.

DE is a population-based evolutionary algorithm (EA), firstly introduced by R. Storn and K. Price in 1997\cite{storn1997differential} as a method to find the global optimum of nonlinear, non-differentiable functions defined over a continuous parameter space.
As indicated by a recent study\cite{das2016recent}, DE exhibits much better performance in comparison with several others continuous-variable meta-heuristic algorithms on a wide range of real-world optimization problems, despite its simplicity.
Being inspired by evolution of species, it exploits the operations of \textit{crossover}, \textit{mutation} and \textit{selection} to generate new candidate solutions. 
However, unlike traditional EAs and GAs, the mutated solutions are generated as scaled differences of distinct individuals of the current population. This self-referential mutation tends to automatically adapt the different variables of the problem to the natural scale of the solution landscape, so improving the search potential of the algorithm. For this reason, it has been elected as optimization engine of the \textit{EOS} algorithm. 

Three major improvements have been made in \textit{EOS} to the standard DE algorithm in order to deal with hard, unconstrained, global optimization problems: (i) a self-adaptation of control parameters, (ii) an epidemic mechanism, to avoid premature convergence in local optima, and (iii) an island-model, to achieve a nice balance between \textit{exploration} and \textit{exploitation} of the search space.

\subsubsection{Self-Adaptation.}
The $j$DE self-adaptive scheme proposed by Brest et al.\cite{brest2007performance} is implemented in \textit{EOS} for automatically adjusting the values of both the scale factor $F$ and the crossover probability $C_r$, that are
the only two control parameters in DE, apart from the population size $N_p$. 
In $j$DE, 
each individual feature its own private copy of the 
control parameters $F$ and $C_r$,
which are randomly initialized and thus different from individual to individual.
Therefore each individual mutates according to its own set of rule. The hope is that ``good'' values of these control parameters would contribute to produce ``better'' individuals, which, being more prone to survive and produce offspring, will propagate their control parameters into the population on successive generations.

\subsubsection{Epidemic Mechanism.}
\begin{figure} [h!]
    \centering
    \includegraphics[width = 0.6\textwidth]{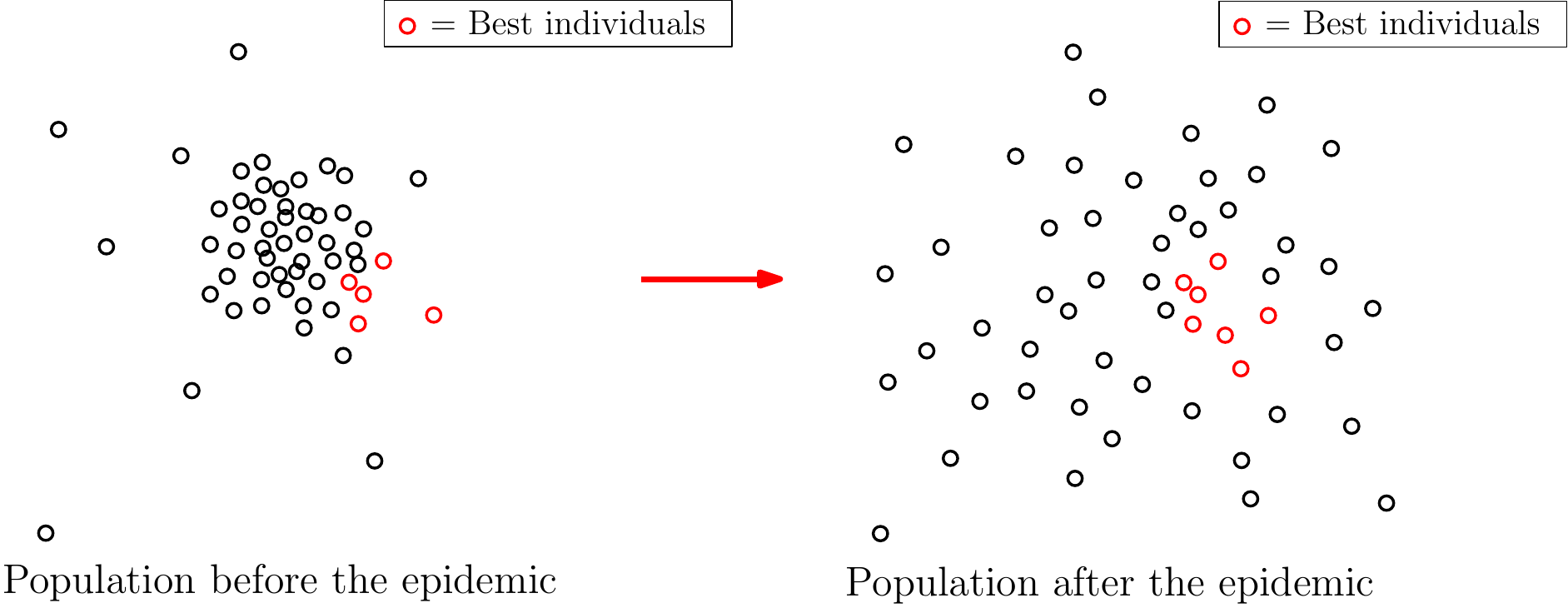}
    \caption{Effect of the Epidemic mechanism on the population.}
    \label{fig:epidemic}
\end{figure}

A partial-restart mechanism, named ``Epidemic'', is adopted in \textit{EOS} so that to maintain diversity between the individuals of the population along generations, that is a fundamental issue in DE. For this propose, population ``diversity'' is evaluated at the end of each generation, in term of Euclidean distance between any pairs of solutions. If the diversity score falls under a certain threshold for the large part of the population and epidemic outbreak. 
A few best individuals are spared, whereas a large portion of the remaining population is randomly reinitialized over the entire search space. This mechanism is illustrated in Figure~\ref{fig:epidemic}.

\begin{figure} [h!]
    \centering
    \includegraphics[width = 0.6\textwidth]{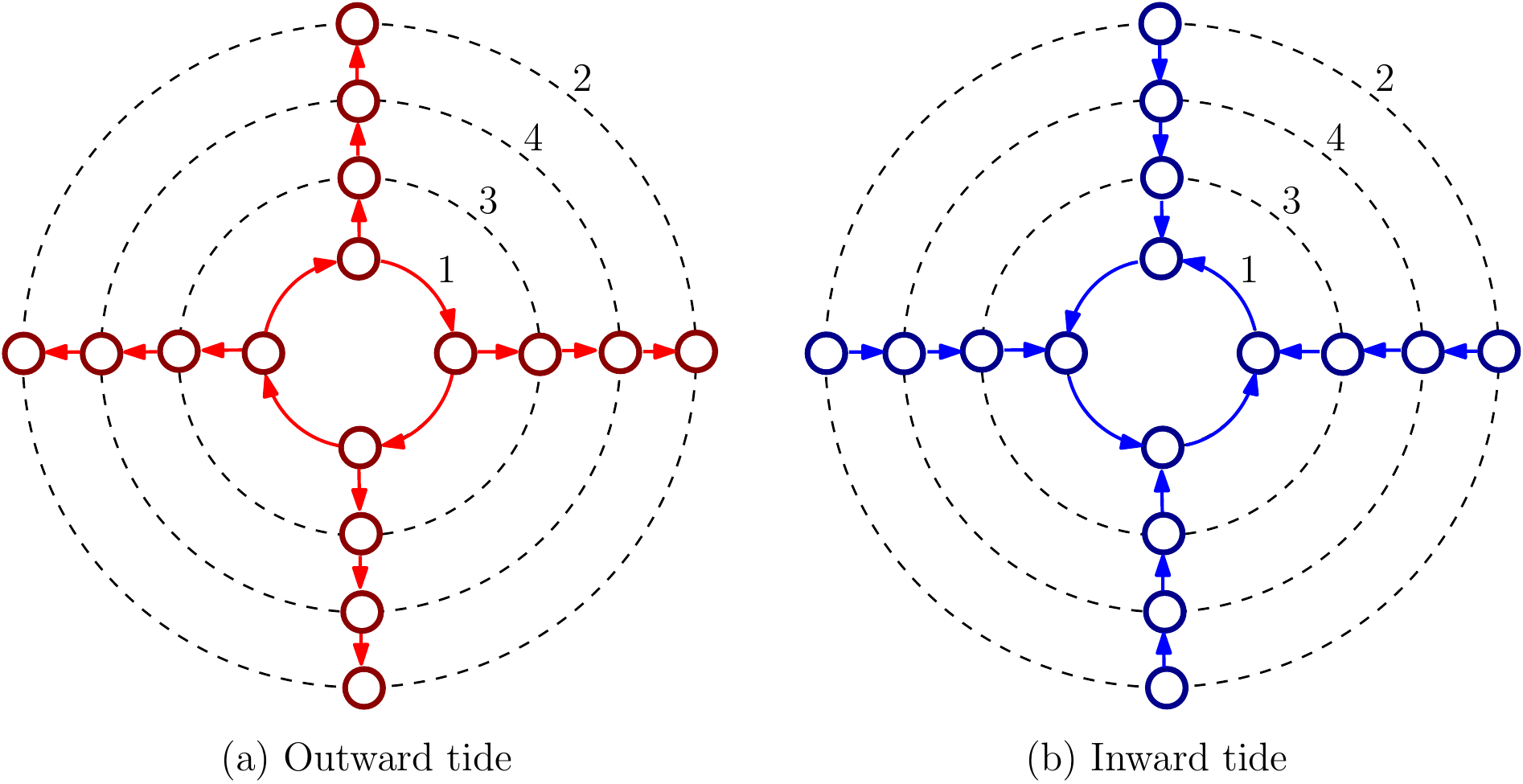}
    \caption{Migration tide: forward (a) and backward (b), for the 16-island case.}
    \label{fig:island}
\end{figure}

\subsubsection{Island Model.}
A number of mutation strategies for DE is documented in literature \cite{das2016recent}, that differ for the number and type of individuals employed in the scaled differences.
In the present application, the following mutation strategies have been implemented: \textit{DE/rand/1/bin} (or Strategy 1), \textit{DE/best/1/bin} (or Strategy 2), \textit{DE/target-to-best/1/bin} (or Strategy 3), \textit{DE/best/2/bin} (or Strategy 4), where the current nomenclature was used\cite{storn1997differential}.
Typically, strategies based on perturbation of the best individual (strategies 2 and 4) show a faster converge rate toward an (often local) minimum (\emph{i.e.}, have a greater \textit{exploitative} power), while strategies based on randomly chosen individuals (strategies 1 and 3) provide a better, yet slower, way to explore the whole search space (\emph{i.e.}, have a greater \textit{explorative} power). 
Some strategies are  better than others on some optimization problems; but
the opposite may be true on different problems.
Therefore,
the need arises to combine the strengths of each of the strategies so that to obtain a more robust and performing algorithm.
The problem has been tackled in \textit{EOS} by exploiting a synchronous island-model: several populations, or ``tribes'', are created, each one evolving on a different ``island'' of an archipelago, with a customized topology.
Each tribe evolves independently according to its own (pre-assigned) mutation rule, until a ``migration'' is performed.
At that point the best $n_b$ agents of each tribe migrate (i.e., are copied)  in the ``connected tribes'', where they replace the $n_b$ worst agents. 
Using heterogeneous strategies along different islands guarantees a good balance between the search space exploration and exploitation, usually overcoming the reliability and performance of the best of the strategies involved \cite{izzo2009parallel}.
In addition, the island-model paradigm allows for an easy parallelization: each island can be assigned to a different thread of a cluster (\emph{e.g.}, through the MPI protocol); the evolution phase proceeds in parallel, until communications between processes are performed during migrations.
The ``migrations'' tides alternate their direction at each event, as shown in Figure~\ref{fig:island} for a radially-arranged 16-island archipelago; in the reported topology, inner rings promote exploration, featuring strategies 1 and 3, while outer rings promote exploitation, featuring strategies 2 and 4.

\section{Simulated Annealing}\label{sez:SA}



Simulated Annealing (SA) is used in the present paper to solve the outer-level combinatorial optimization problems.
SA is a simple, yet powerful, neighborhood-based meta-heuristic algorithm for global optimization. 
Inspired by a strong analogy between the search for the global optimum of a black-box function and the
search for the minimum-energy configuration of a solid material, that is, a ``physical'' annealing process,
SA has been massively applied to real-life applications,
especially in case the
evaluation of the objective function requires a large amount of memory
or to deal with combinatorial optimization problems,
where gradient-like information provided by population-based algorithms are usually less relevant.
A brief overview of Simulated Annealing is here presented. Interested readers are suggested to refer to Reference~\citenum{Aarts} for further details.
%
%
A pseudo-code of the basic algorithm is proposed in Figure~\ref{alg:SA-alg}, highlighting the presence of three fundamental operators, namely neighbor generation function, acceptance rule, and cooling schedule, which form the backbone of any SA.

At the beginning an initial solution $x$, often referred as \textit{system configuration} or  \textit{state},  is randomly generated.
Its fitness $f(x)$, or \textit{energy},  is evaluated according to the objective function to minimize.
Then,
a new solution $x' = g_N(x)$ 
is generated in a neighborhood $N$ of $x$, that is, $x'$ is a perturbation of the current solution $x$ according to some ``neighbor'' generation function $g_N$.
The new solution is then compared to the previous one in order to decide whether the transition to the new configuration should be accepted.
The Metropolis criterion,
originally devised for describing the state transition probability in a physical annealing process, 
is almost-unanimously adopted as acceptance rule in SA.
According to this criterion, a transition is always accepted if it reduces the energy of the system, that is, if $f(x') < f(x)$.
Otherwise, $x'$ could also be accepted with probability $p=exp\left(-\frac{f(x') - f(x)}{T}\right)$, where $T$ is a control parameter, or \textit{temperature}, which 
allows for tuning the acceptance ratio over the course of the search process.

The temperature is updated at any iteration according to a prescribed law, or \textit{cooling schedule}, which brings this control parameter from an initial value $T_0$ to a final one $T_F$ in a certain number of steps.
At the beginning of the search, a sufficiently high temperature value is needed in order to allow the exploration of the whole search space and (potentially) reach solutions not in the immediate neighborhood of the initial  system configuration.
In the latter stages of the search process, lower values are needed, limiting the explorative capability of the algorithm while favouring local improvements.
This process is repeated iteratively until some termination criterion,
that is, a maximum number of iterations or a very small temperature, is met.

\begin{center}
\begin{minipage}{.6\linewidth}
\centering
\begin{algorithm} [H] \small
        \caption{Simulated Annealing}
        \begin{algorithmic}
        \Procedure{SA}{$f$}
 	\State $x = x_0$ \Comment Initialize solution
 	\State $T = T_0$ \Comment Initialize temperature
 	\While{ $T < T_F$ } 
 	    \State $x' \gets g_N(x)$  \Comment \textbf{Neighbor generation}
 		\State $\Delta \gets f(x')-f(x)$ 
 		\If {$\Delta \le 0$} \Comment \textbf{Acceptance rule}
            \State $x\gets x'$
        \Else
           \State$x \gets x'$ with probability $e^{-\Delta/T}$
        \EndIf
        \State $T\gets updateTemperature()$ \Comment \textbf{Cooling schedule}
	\EndWhile
	\State \textbf{return $x$}
	\EndProcedure
\end{algorithmic}
\end{algorithm}
\captionof{figure}{Simulated annealing basic algorithm}
\label{alg:SA-alg}
\end{minipage}
\end{center}


\subsection{Cooling Schedule}\label{sez:cooling}
A proper cooling schedule is of paramount importance for attaining good-quality solutions.
Under mild assumptions, asymptotic convergence towards the global optimum can be proved for certain, very slow, cooling schedules, such as the \textit{Logarithmic Schedule}:
$T = b/log(a+k)$, with $a$, $b$ positive real numbers and $k$ the iteration number.
Unfortunately, the computational time required by this scheduling
to reach convergence (i.e., a sufficiently low temperature) is very high.
For this reason, several faster cooling schedules have been proposed in literature on the basis of practitioners' experience with typical problems.
Among them, the most common choice is given by the \textit{Geometrical Schedule}: 
\begin{equation}
\label{eq:geoCool}
    T = T_0 \, \alpha^{\left\lfloor  k/L\right\rfloor}  
\end{equation} 
where
$k$ is the iteration number, 
$\alpha$ the cooling factor and 
$L$ the number of steps at a constant temperature.

A proper tuning of the problem-dependent scheduling parameters ($T_0$, $\alpha$, $L$) is usually needed to attain good results. 
Theoretical analysis suggest that $L$ should be of the same order as the neighbor size (hence it depends on the problem and the generating function).
The other parameters requires a more laborious trial and error procedure, with parameters typically in the ranges
$T_0 \in [0.5, 20]$, 
$\alpha \in [0.7, 0.999]$,
and $T_F \in [10^{-7}, 10^{-5}]$.

%

%



\subsection{Neighbor generation function}\label{sez:cooling}
\begin{figure} [h!]
    \centering
    \subfigure[Insert.]{ \label{fig:insert} \includegraphics[width=0.25\textwidth]{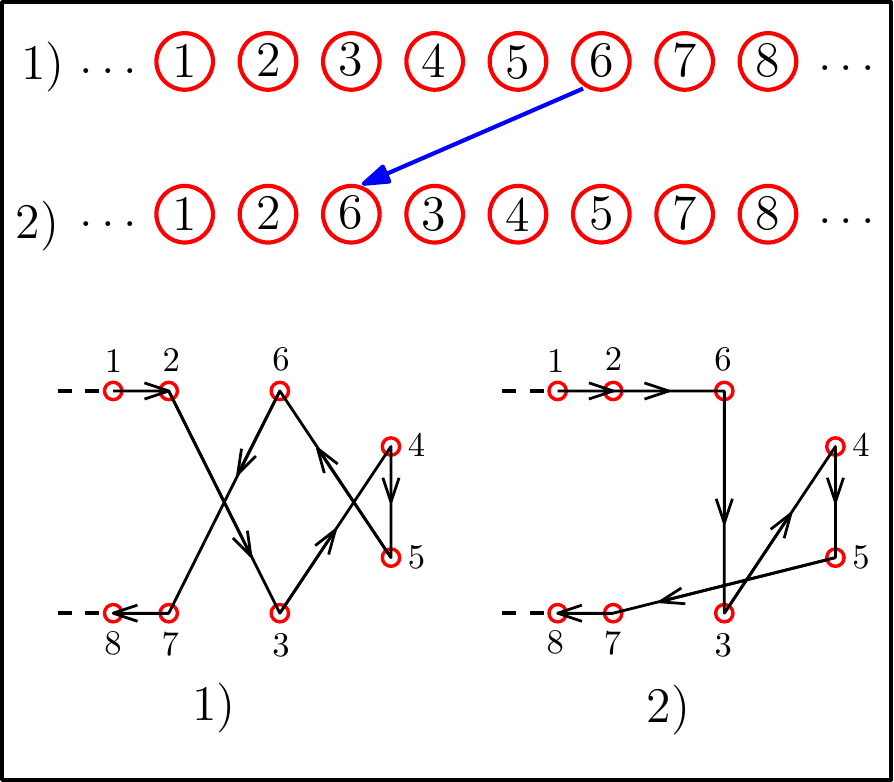}}
    $\qquad$
    \subfigure[Reverse.]{ \label{fig:reverse} \includegraphics[width=0.25\textwidth]{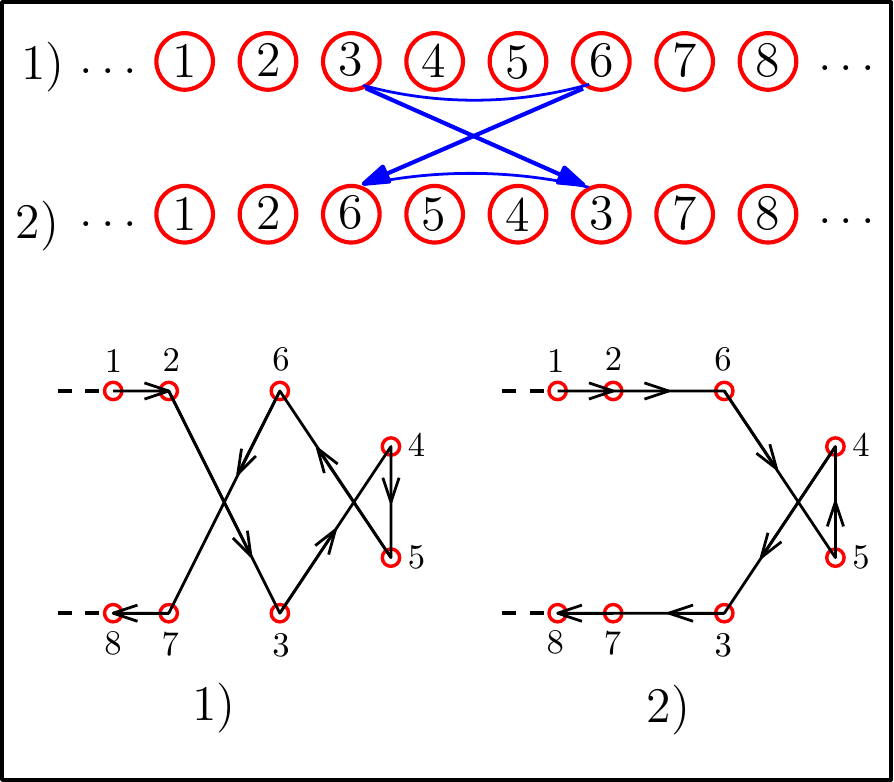}}\\
    \subfigure[Swap.]{ \label{fig:swap} \includegraphics[width=0.25\textwidth]{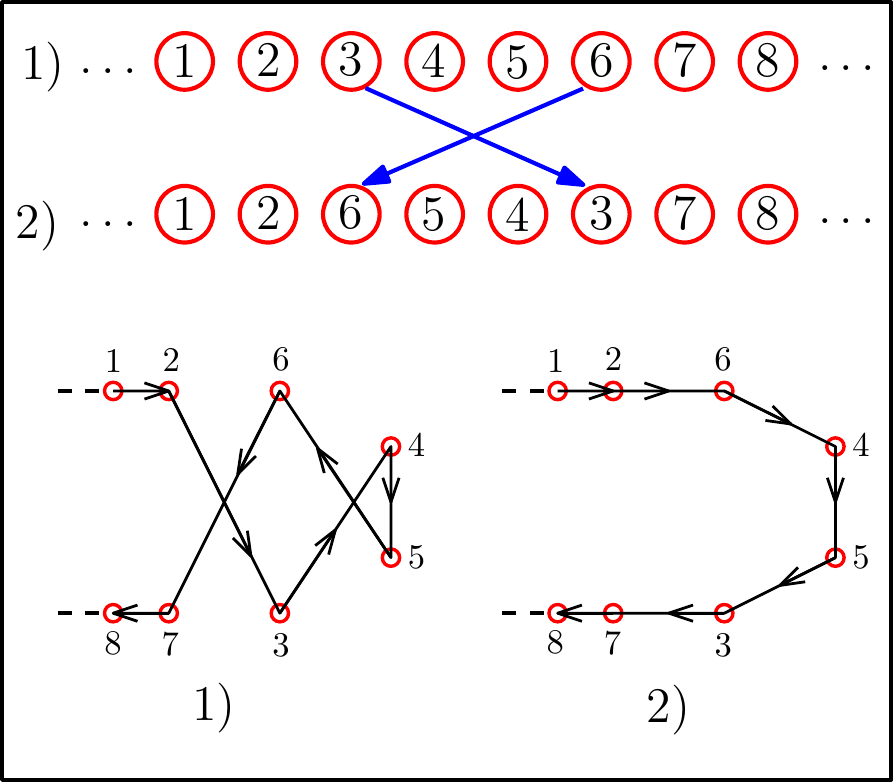}}
    $\qquad$
    \subfigure[Scramble.]{ \label{fig:scramble} \includegraphics[width=0.25\textwidth]{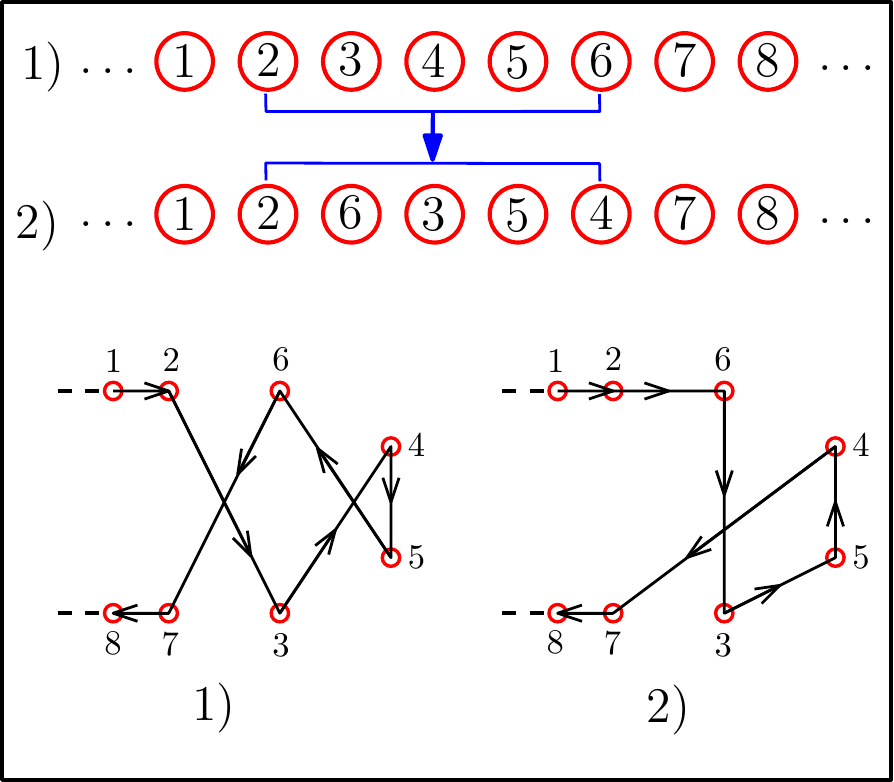}}
    \caption{Representation of the neighbor generation functions in terms of the effect on the permutation-encoded state $x$ and on the corresponding TSP path.}
    \label{fig:Moves}
\end{figure}

The neighbor generation is tightly related to the adopted \textit{encoding}, that is, the way the real-word problem is described in terms of numerical variables.
If a permutation encoding is considered, a local perturbation of the current solution
can be obtained by using permutation-preserving mutation operators developed for solving the classic TSP problem,\cite{Kirkpatrick1983}
such as i) insert, ii) swap, iii) reverse, and iv) scramble. 
A visual representation in terms of the corresponding TSP path change is proposed in  Figure~\ref{fig:Moves}.
Each operator defines a different neighborhood of the current solution $x$, that is,
a set $N(x)$ composed by all the immediate ``neighbors'' that can be generated from $x$.
Consequently, each operator may be more or less effective on a given problem.






\section{Numerical Results}

In this section numerical results are proposed for various touring missions of pre-determinate length.
Orbital parameters of both chaser and targets are provided in Table~\ref{tab:targets}.
\begin{table}[!h]
    \centering
    \resizebox{0.8\textwidth}{!}{   
    \begin{tabular}{l c | c c c c c c c c c c} 
\hline 
 & \textbf{Chaser} & \multicolumn{10}{c}{\textbf{Targets}} \\ 
\hline 
ID & 0 & 1 & 2 & 3 & 4 & 5 & 6 & 7 & 8 & 9 & 10  \\ 
$r\, [\SI{}{\kilo\m}]$ & 7000 & 6900 & 6910 & 6930 & 6940 & 6950 & 6960 & 6980 & 7010 & 7020 & 7030 \\ $\theta_0\, [deg]$ &  0 & -5 & 10 & 15 & 35 & -30 & -10 & 25 & 20 & -25 & -15\\ 
\hline
ID &  & 11 & 12 & 13 & 14 & 15 & 16 & 17 & 18 & 19 & 20  \\ 
$r\, [\SI{}{\kilo\m}]$ &   & 7050& 7060& 7070& 7080&	7090&	7110&	7130&	7140&	7160&	7170 \\
$\theta_0\, [deg]$ &    & 5&	-35&	30&	-20&	-40&	50&	-50&	40&	-45&	45\\ 

\hline 
\end{tabular} 

    }
    \caption{Chaser and targets initial orbital parameters.}
    \label{tab:targets}
\end{table}

\begin{figure}[!h]
 	\centering
 	\subfigure[The optimal value and the heuristic estimate of the $\Delta V$ for a transfer leg from $r_1=$ \SI{7000}{\kilo\metre} towards $r_2=$ \SI{7140}{\kilo\metre} is reported as a function of the phase angle $\Delta \theta$ at the departure and the transfer time $t$, expressed as a multiple of the departure orbit period $T_0$. The $\Delta V$ of the Hohmann transfer is also reported.]
 	{\label{fig:heuComp-DV}
 		 \includegraphics[width=0.7\textwidth]{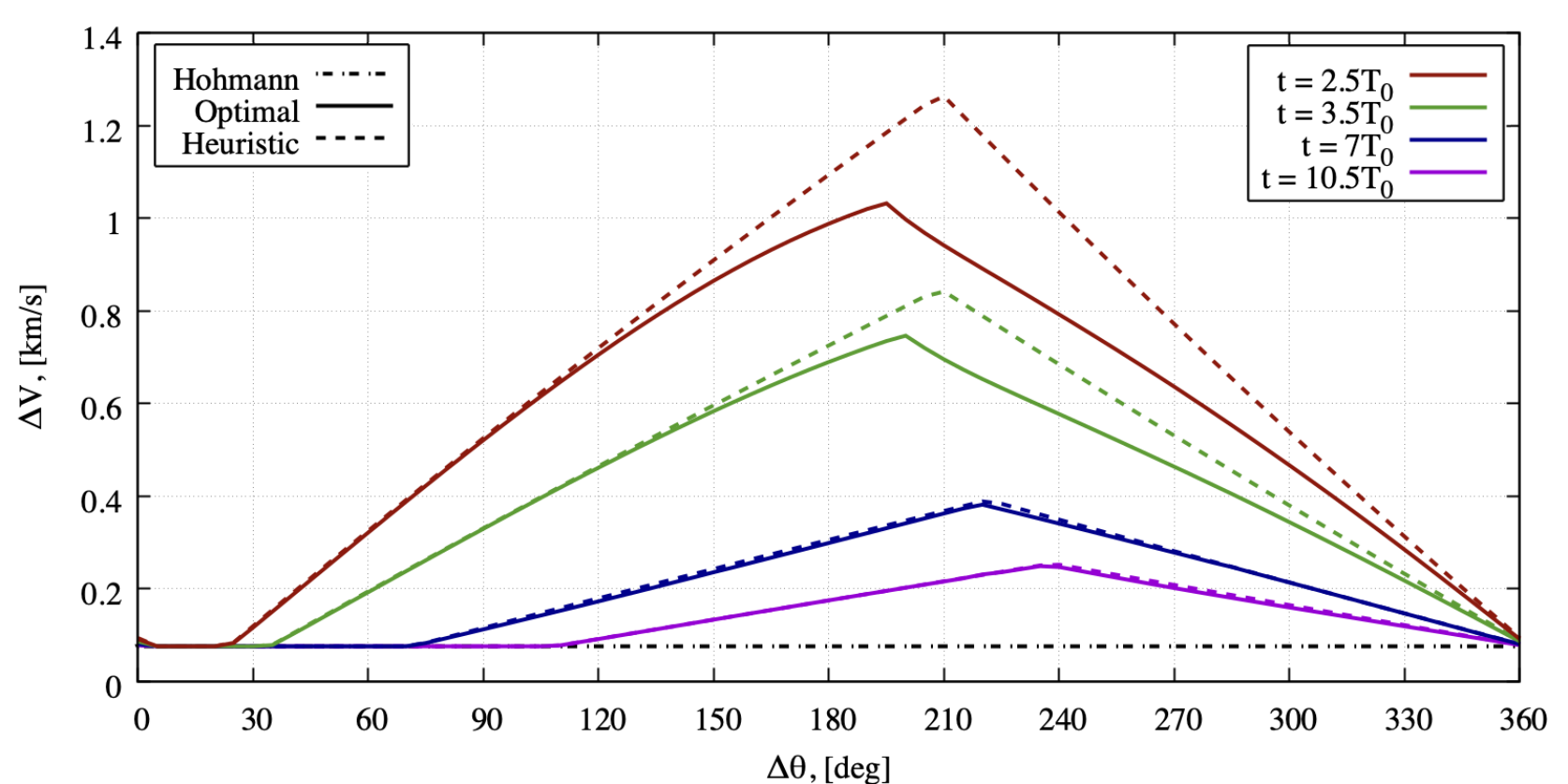}
 		} \\
 	\subfigure[Optimal and sub-optimal chaser trajectory for $t = 2.5T_0$ and $\Delta \theta_0 = 150 \degree$]
 	{\label{fig:heuComp-traj}
 		    \includegraphics[width = 0.7\textwidth]{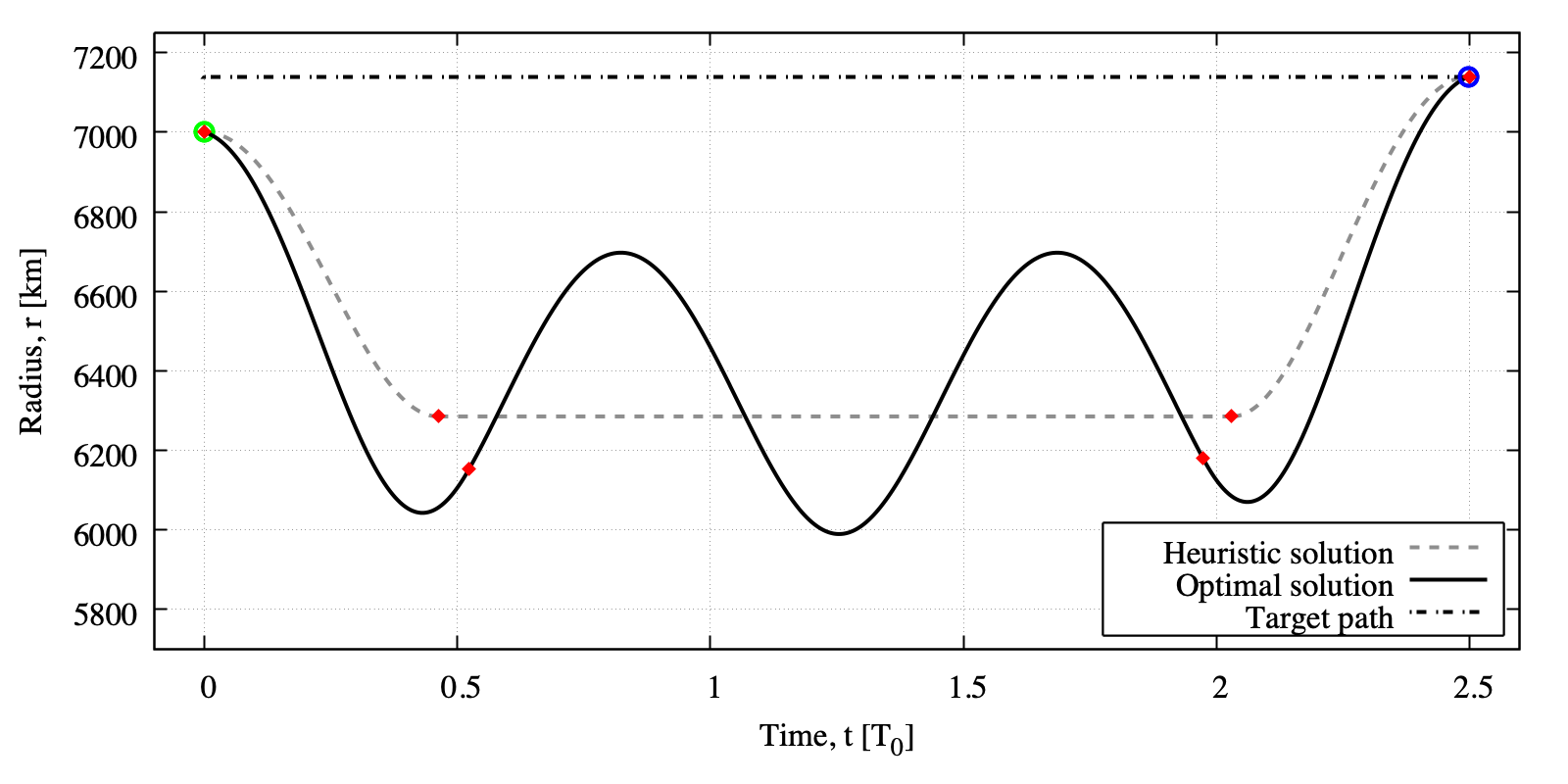}
    }
 	\caption{Comparison between the optimal chaser trajectory and the sub-optimal mission scenario elected as heuristic.}
 	\label{fig:heuComp}
\end{figure}

%
The validity of the sub-optimal, analytic, four-impulse, solution for the co-planar rendezvous has been investigated first,
by comparing it with the solution obtained after a full optimization procedure, in a few relevant cases.
Figure~\ref{fig:heuComp-DV} shows the optimal and (heuristic) estimated $\Delta V$ for a transfer leg with
departure orbit radius $r_1=$ \SI{7000}{\kilo\metre} and 
arrival orbit radius $r_2=$ \SI{7140}{\kilo\metre},
as a function of the phase angle $\Delta \theta$ at the departure,
for transfer times equal to $2.5$, $3.5$, $7$, and $10$ times the departure orbit period $T_0$.
It is apparent that the heuristic solution provides an estimate very close to the optimal value;
typical differences between the optimal and the heuristic $\Delta V$
are in the range of $0.2 \div  3\%$, with few peaks at $20\%$ when the travel time becomes too short for the proposed mission scheme to be practical.
The absolute values of the difference usually do not exceed a few tens of \SI{}{\meter/\second}, and thus are deemed reasonable.

Differences between the optimal and sub-optimal solutions
are apparent in the spacecraft trajectory,
when an unfavourable phasing condition is considered, as shown in Figure~\ref{fig:heuComp-traj}.
Indeed, in contrast with the sub-optimal scheme, the optimal solution adopts an elliptical waiting orbit in order to reduce the cost of the internal burns. Luckily, major discrepancies arise only for very unfavourable departure phases, when even the cost of the optimal maneuver is very high and it is unlikely to show up in the optimal multi-rendezvous solution.

%


\begin{figure} [!h]
    \centering
    \includegraphics[width = 0.7\textwidth]{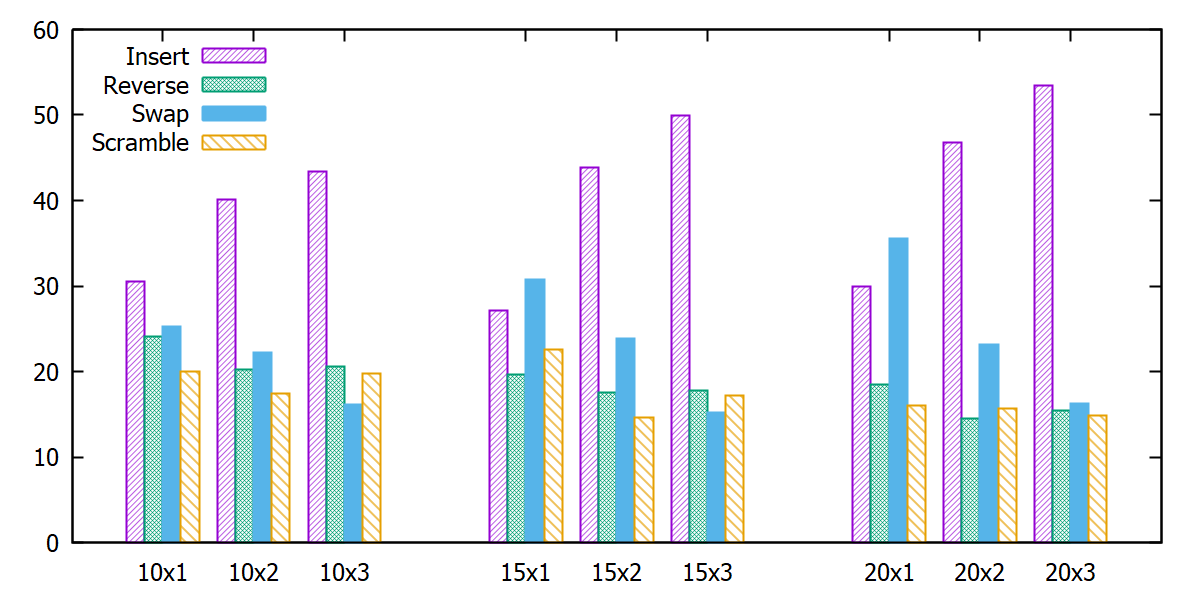}
    \caption{Relative frequency of improving solutions produced by each operator, averaged over  25 independent runs, for the 10-, 15-, and 20-target missions.}
 	\label{fig:SA-moves-res}
\end{figure}

Figure~\ref{fig:SA-moves-res} shows the relative frequency of improving solutions produced by each operator, averaged over  25 independent runs, for different missions.
For the $N \times 1$ tours ($N \in \{10, 15, 20\})$, all entries in the permutation correspond to a valid target, and all operators seems to perform equally well.
Instead, when the number of divisions is larger, 
most of the entries in a permutation are blanks; as a result,
the \textit{insert}  operator becomes more efficient and 
more frequently leads to an improving solution.

A comparison between the minimum $\Delta V$ obtained
at different stage of the optimization process (outer-level, inner-level with fixed or free encounter times) is presented in Tables~\ref{tab:comp}(a) and \ref{tab:comp}(b) for the 10-target and 20-target missions, respectively.
Details on the best attained solutions for these missions are reported in Tables~\ref{tab:Sol8} and
\ref{tab:Sol15}. 
Differences in total $\Delta V$ estimated by SA and fixed-times DE are in the order of few percents, confirming the preliminary analysis on the effectiveness of the employed heuristic.
The saving is far greater (up to $20\%$) if the DE algorithm is left free to modify the encounter times, as the algorithm is able to perform a more appropriate sizing of the time window allocated for each leg, shortening or extending legs as needed. As expected, these differences decreases  
when the number of time divisions $D$ is larger.
As a final remark, it is worth to note that in the problem under investigation 
the mission time is always proportional to the number of targets, 
hence the latter optimization step is more likely to be  important for small size problem, when the overall mission duration is lower and a smaller number of  favourable (i.e., Hohmann) transfer opportunities exists.

\begin{figure} [h!]
    \centering
    \includegraphics[width = 0.8\textwidth]{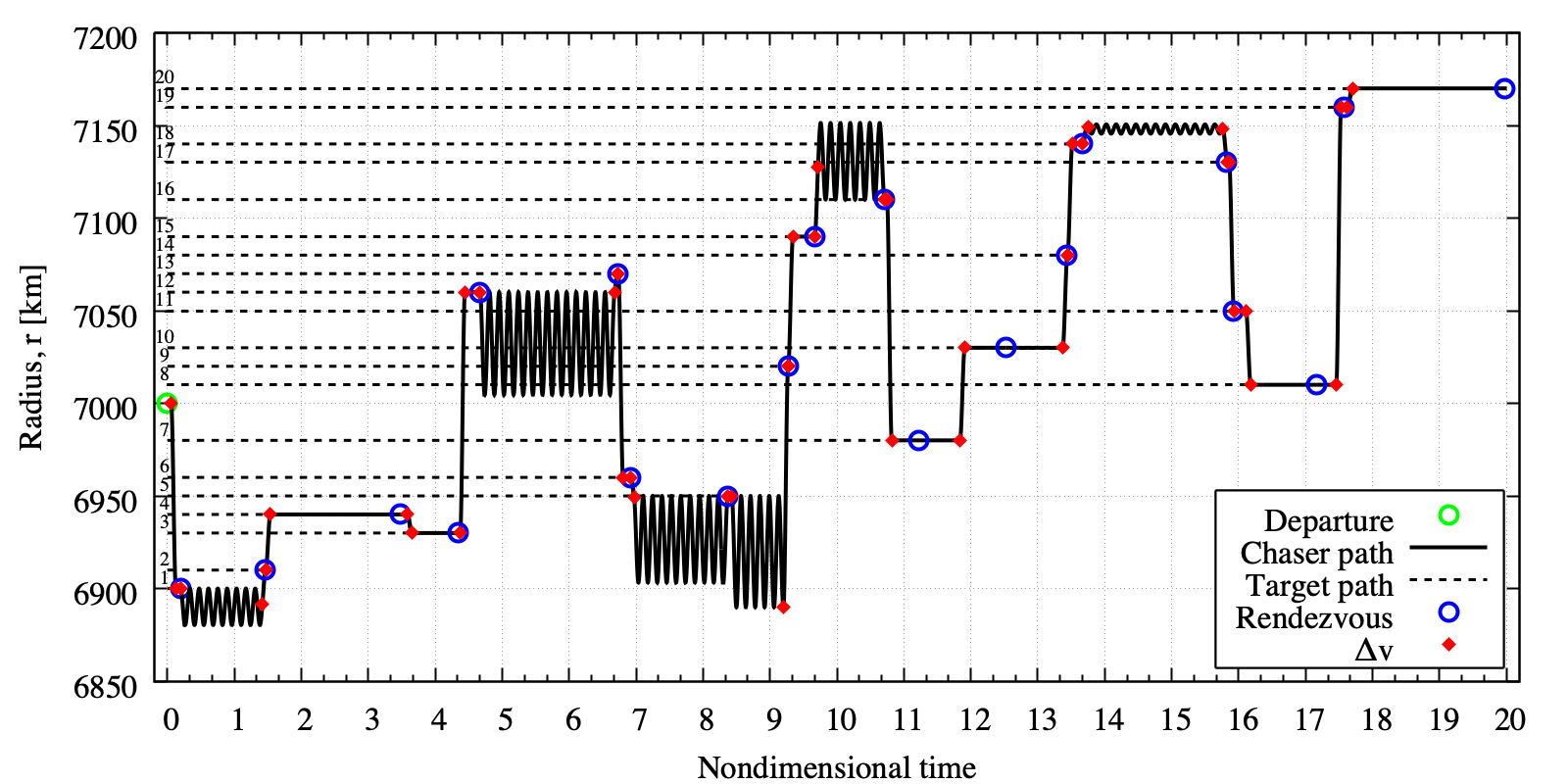}
    \caption{Radius $r$ vs. time $t$ for the 20x3 solution.}
    \label{fig:Rt15x3}
\end{figure}

\begin{table}[h!]
 \centering
    \subtable[10-target missions]{
    \resizebox{0.4\textwidth}{!}{
    \begin{tabular}{l c c c } 
\hline 
\hline 
& \textbf{SA} & \textbf{DE} & \textbf{DE} \\ 
& & time-fixed & time-free \\ 
\hline 
\textbf{10x1}: $\Delta V_{tot}\, [\SI{}{\kilo\m}/\SI{}{\s}]$ & 0.6181 & 0.5948 & 0.4980 \\ 
\textbf{10x2}: $\Delta V_{tot}\, [\SI{}{\kilo\m}/\SI{}{\s}]$ & 0.4828 & 0.4788 & 0.4726 \\ 
\textbf{10x3}: $\Delta V_{tot}\, [\SI{}{\kilo\m}/\SI{}{\s}]$ & 0.4698 & 0.4662 & 0.4488 \\ 
\hline 
\hline 
\end{tabular} 
}
    \label{tab:8tarcomp} }
    \subtable[20-target missions]{
    \resizebox{0.4\textwidth}{!}{ 
    \begin{tabular}{l c c c } 
\hline
\hline 
& \textbf{SA} & \textbf{DE} & \textbf{DE} \\ 
& & time-fixed & time-free \\ 
\hline 
\textbf{20x1}: $\Delta V_{tot}\, [\SI{}{\kilo\m}/\SI{}{\s}]$ & 0.8815 & 0.8688 & 0.7592 \\ 
\textbf{20x2}: $\Delta V_{tot}\, [\SI{}{\kilo\m}/\SI{}{\s}]$ & 0.7899 & 0.7819 & 0.7506 \\ 
\textbf{20x3}: $\Delta V_{tot}\, [\SI{}{\kilo\m}/\SI{}{\s}]$ & 0.7715 & 0.7651 & 0.7449 \\ 
\hline 
\hline
\end{tabular} 
}
    \label{tab:15tarcomp}}
    \caption{Attained solutions at various stages of the optimization procedure.
    SA refers to the solution of the outer-level problem, DE (time-fixed) and DE (time-free)
    refer to the refined solution attained after completing the inner-level optimization, assuming the encounter epochs respectively fixed or free to be optimized.}
    \label{tab:comp}
\end{table}

\begin{table}[h!]
    \caption{Optimal 10-target sequences.}
    \centering
    \resizebox{0.7\textwidth}{!}{
    \begin{tabular}{l c c c c c c c c c c } 
\hline 
\hline 
\multicolumn{11}{c}{\textbf{Mission 10x1:} $\Delta V_{tot} = 0.49796\, \SI{}{\kilo\m}/\SI{}{\s}$} \\ 
\hline 
ID & 8 & 7 & 1 & 2 & 3 & 4 & 9 & 10 & 5 & 6 \\ 
$t\, [\SI{}{\day}]$ & 0.2419 & 0.6163 & 1.1430 & 1.5159 & 2.3611 & 2.9687 & 3.2728 & 3.8538 & 4.1485 & 4.7222 \\ 
$\Delta V\, [\SI{}{\kilo\m}/\SI{}{\s}]$ & 0.0836 & 0.0557 & 0.0950 & 0.0055 & 0.0468 & 0.0176 & 0.0433 & 0.0525 & 0.0432 & 0.0547 \\ 
\hline 
\hline 
\multicolumn{11}{c}{\textbf{Mission 10x2:} $\Delta V_{tot} = 0.47261\, \SI{}{\kilo\m}/\SI{}{\s}$} \\ 
\hline 
ID & 6 & 8 & 7 & 1 & 2 & 3 & 4 & 10 & 9 & 5 \\ 
$t\, [\SI{}{\day}]$ & 0.2338 & 0.4957 & 0.7821 & 1.3828 & 1.6548 & 2.2590 & 2.8333 & 3.0230 & 4.4469 & 4.7222 \\ 
$\Delta V\, [\SI{}{\kilo\m}/\SI{}{\s}]$ & 0.0217 & 0.0374 & 0.0972 & 0.1145 & 0.0109 & 0.0640 & 0.0189 & 0.0487 & 0.0217 & 0.0379 \\ 
\hline 
\hline 
\multicolumn{11}{c}{\textbf{Mission 10x3:} $\Delta V_{tot} = 0.44876\, \SI{}{\kilo\m}/\SI{}{\s}$} \\ 
\hline 
ID & 6 & 8 & 5 & 3 & 4 & 1 & 9 & 10 & 7 & 2 \\ 
$t\, [\SI{}{\day}]$ & 0.2902 & 0.6824 & 0.7402 & 1.6528 & 1.7315 & 2.3959 & 2.5057 & 3.1332 & 4.5675 & 4.7222 \\ 
$\Delta V\, [\SI{}{\kilo\m}/\SI{}{\s}]$ & 0.0217 & 0.0270 & 0.0324 & 0.0761 & 0.0054 & 0.0818 & 0.0653 & 0.0340 & 0.0668 & 0.0382 \\ 
\hline 
\hline 
\end{tabular} 

    }
    \label{tab:Sol8}
\end{table}

\begin{table}[h!]
    \caption{Optimal 20-target sequences.}
    \centering
    \resizebox{\textwidth}{!}{
    \begin{tabular}{l c c c c c c c c c c c c c c c c c c c c } 
\hline
\hline 
\multicolumn{21}{c}{\textbf{Mission 20x1:} $\Delta V_{tot} = 0.75921\, \SI{}{\kilo\m}/\SI{}{\s}$} \\ 
\hline 
ID & 1 & 2 & 4 & 3 & 12 & 13 & 6 & 5 & 9 & 16 & 10 & 7 & 14 & 18 & 17 & 11 & 8 & 20 & 19 & 15 \\ 
$t\, [\SI{}{\day}]$ & 0.2847 & 0.6928 & 1.6389 & 1.9387 & 2.1517 & 2.8333 & 3.2135 & 3.8157 & 4.1377 & 4.4745 & 4.8359 & 5.6666 & 6.0980 & 6.4059 & 7.3194 & 7.5257 & 8.1069 & 8.3114 & 9.2043 & 9.4444 \\ 
$\Delta V\, [\SI{}{\kilo\m}/\SI{}{\s}]$ & 0.0545 & 0.0218 & 0.0165 & 0.0054 & 0.0702 & 0.0522 & 0.0591 & 0.0327 & 0.0379 & 0.0478 & 0.0425 & 0.0269 & 0.0536 & 0.0316 & 0.0161 & 0.0423 & 0.0214 & 0.0846 & 0.0052 & 0.0367 \\ 
\hline 
\hline 
\multicolumn{21}{c}{\textbf{Mission 20x2:} $\Delta V_{tot} = 0.75083\, \SI{}{\kilo\m}/\SI{}{\s}$} \\ 
\hline 
ID & 1 & 2 & 4 & 3 & 12 & 13 & 6 & 5 & 9 & 16 & 15 & 10 & 7 & 14 & 18 & 17 & 11 & 8 & 19 & 20 \\ 
$t\, [\SI{}{\day}]$ & 0.1181 & 0.6928 & 1.5814 & 1.7440 & 2.1250 & 3.1631 & 3.2359 & 3.6992 & 4.3083 & 4.4861 & 4.9707 & 5.5170 & 5.7582 & 6.1435 & 6.3749 & 7.4374 & 7.5471 & 8.0364 & 8.2963 & 9.4444 \\ 
$\Delta V\, [\SI{}{\kilo\m}/\SI{}{\s}]$ & 0.0545 & 0.0175 & 0.0165 & 0.0054 & 0.0702 & 0.0331 & 0.0591 & 0.0415 & 0.0682 & 0.0478 & 0.0291 & 0.0319 & 0.0269 & 0.0536 & 0.0336 & 0.0136 & 0.0423 & 0.0214 & 0.0794 & 0.0052 \\ 
\hline 
\hline 
\multicolumn{21}{c}{\textbf{Mission 20x3:} $\Delta V_{tot} = 0.74495\, \SI{}{\kilo\m}/\SI{}{\s}$} \\ 
\hline 
ID & 1 & 2 & 4 & 3 & 12 & 13 & 6 & 5 & 9 & 15 & 16 & 7 & 10 & 14 & 18 & 17 & 11 & 8 & 19 & 20 \\ 
$t\, [\SI{}{\day}]$ & 0.0921 & 0.6928 & 1.6452 & 2.0481 & 2.2037 & 3.1793 & 3.2681 & 3.9528 & 4.3780 & 4.5648 & 5.0580 & 5.3073 & 5.9131 & 6.3476 & 6.4536 & 7.4768 & 7.5271 & 8.1048 & 8.2953 & 9.4444 \\ 
$\Delta V\, [\SI{}{\kilo\m}/\SI{}{\s}]$ & 0.0545 & 0.0164 & 0.0165 & 0.0054 & 0.0702 & 0.0349 & 0.0591 & 0.0310 & 0.0708 & 0.0373 & 0.0322 & 0.0694 & 0.0269 & 0.0267 & 0.0316 & 0.0137 & 0.0423 & 0.0214 & 0.0794 & 0.0052 \\ 
\hline 
\hline 
\end{tabular} 

    }
    \label{tab:Sol15}
\end{table}


\section{Conclusions}

This paper investigated a bi-level optimization procedure for the design of the multi-rendezvous trajectory of a chaser spacecraft which has to visit a prescribed set of space debris. The goal was to minimize the overall propellant consumption, while completing the tour within a given amount of time. 
First, analogies between the time-dependent TSP and the combinatorial features of an ADR mission,
i.e., the definition of the optimal encounter sequence together with a preliminary evaluation of the epochs at each encounter, have been highlighted.
Two alternative formulations, the Integer Linear Programming and the Permutation Optimization, are recalled for the TSP and extended in order to deal with cost functions depending on both departure time and transfer duration.
Simulated Annealing was used to solve the permutation optimization problems previously defined. A simple, sub-optimal, analytic solution of the single-target rendezvous problem was adopted as heuristic for a fast evaluation of the $\Delta V$ associated to each leg, without studying it in full details.
The attained solutions were further refined by assuming the encounter sequence fixed and optimizing the multi-impulse rendezvous trajectory: each body-to-body transfer was described by means of a peculiar parameterization based on the position of the impulses, whose magnitude was minimized.
Numerical solutions were presented for a set comprising up to 20 target bodies. Results suggested that, by coupling the proposed $\Delta V$ heuristic with a time-discrete time-fixed formulation with time-discretization factor $3$, one attains a trajectory that is very close to the solution of the full mixed-integer nonlinear programming problem, whereas the overall computational effort is significantly reduced.


\bibliographystyle{AAS_publication}   
\bibliography{MRRreferences}   

\end{document}